\documentclass[12pt]{amsart}
\usepackage{amsmath,amsthm,latexsym,amscd,amsbsy,amssymb}
 \relax

\chardef\bslash=`\\ 

\makeatletter
\def\verbatim{\interlinepenalty\@M \@verbatim
 \leftskip\@totalleftmargin\advance\leftskip2pc
 \frenchspacing\@vobeyspaces \@xverbatim}
\makeatother
\makeatletter
 \def\dgt@k{\dg@DX=-3 \dg@DY=2 \dg@SIZE=3}
\makeatother

\makeatletter
 \def\dgt@kk{\dg@DX=3 \dg@DY=-1 \dg@SIZE=3}%
\makeatother

\theoremstyle{plain}
\newtheorem{thm}{Theorem}[section]
\newtheorem{cor}[thm]{Corollary}
\newtheorem{lem}[thm]{Lemma}
\newtheorem{pro}[thm]{Proposition}

\newtheorem*{A}{Theorem 7.6}
\newtheorem*{B}{Theorem 6.4}

\theoremstyle{definition}
\newtheorem{rem}[thm]{Remark}
\newtheorem{defin}[thm]{Definition}

\newcommand{\ed}{\operatorname{edd}}

\numberwithin{equation}{section}

\def\cal{\bf}
\def\calibr{{\rm cal}}

\def\cone{{\rm Cone}}

\def\cov{{\rm cov}}
\def\diam{{\rm diam}}
\def\dim{{\rm dim}}
\def\dist{{\rm dist}}
\def\ed{\text{{\rm e-dim}}}
\def\ep{\varepsilon}

\def\Int{{\rm Int}}
\def\ord{{\rm ord}}
\def\pr{{\rm pr}}
\def\R{{\mathbb R}}
\def\st{{\rm st}}
\def\St{{\rm St}}
\def\wh{\widehat}
\def\wt{\widetilde}
\newcommand{\Sgm}[1]{\Sigma^{(#1)}}

\def\sup{{\rm sup}}
\def\inf{{\rm inf}}

\begin{document}


\title[Approximations and selections of multivalued mappings of
finite-dimensional spaces]
{Approximations and selections of multivalued mappings of
finite-dimensional spaces}
\author{N. Brodsky}
\address{Department of Mathematics and Statistics,
University of Saskatche\-wan,
McLean Hall, 106 Wiggins Road, Saskatoon, SK, S7N 5E6,
Canada}
\email{brodsky@math.usask.ca}
\author{Alex Chigogidze}
\address{Department of Mathematics and Statistics,
University of Saskatche\-wan,
McLean Hall, 106 Wiggins Road, Saskatoon, SK, S7N 5E6,
Canada}
\email{chigogid@math.usask.ca}
\thanks{Second author was partially supported by NSERC research grant.}
\author{A. Karasev}
\address{Department of Mathematics and Statistics,
University of Saskatche\-wan,
McLean Hall, 106 Wiggins Road, Saskatoon, SK, S7N 5E6,
Canada}
\email{karasev@math.usask.ca}
\keywords{Continuous selection; approximation; extension dimension}
\subjclass{Primary: 54C65; Secondary: 54C60, 55M10}


\begin{abstract}{We prove extension-dimensional versions of finite
dimensional selection and approximation theorems.
As applications, we obtain several results on extension dimension.}
\end{abstract}

\maketitle
\markboth{N.~Brodsky, A.~Chigogidze, A.~Karasev}{Approximations
and selections of multivalued mappings of finite-dimensional spaces}

\section{Introduction and Preliminary Definitions}\label{S:intro}

Finite-dimensional selection theorem of E.~Michael is very useful
in geometric topology and it is one of central theorems in the
theory of continuous selections of multivalued mappings~\cite{RS}.
A stronger selection theorem is proved in~\cite{SB} and a technique
of its proof shows an interesting interference between selections
and approximations of multivalued mappings.
In particular, finite dimensional approximation theorem
was used in the proof of selection theorem.
However, approximation theorem itself is widely applicable
in mathematics, not only in topology (see a survey~\cite{Kr1}).

There is a new approach in dimension theory exploiting a notion
of extension dimension~\cite{Dr2},\cite{DrD}.
Let $L$ be a CW-complex.
A space $X$ is said to have {\it extension dimension} $\le [L]$
(notation: $\ed X\le [L]$) if
any mapping of its closed subspace
$A\subset X$ into $L$ admits an extension to the whole space
$X$\footnote{Everywhere below $[L]$ denotes the class of
complexes generated by $L$ with respect to the above extension
property, see~\cite{Dr2}, \cite{DrD}, \cite{Ch} for details.}. It
is clear that $\dim X\leq n$ is equivalent to
$\ed X\leq\left[ S^{n}\right]$.

The main purpose of this paper is to prove an extension-dimensional
versions of finite dimensional selection and approximation theorems.
Of course, these versions have the original finite dimensional
theorems as a partial cases. And our proofs follow the ideas
from the paper~\cite{SB}. There is an extension dimensional
approximation theorem for mappings of C-space~\cite{BCh}. We are
mainly interested in the separable
and metrizable situation. In the meantime proofs of our statements without
significant complications remain valid in a more general case of
paracompact spaces and we state our results for the latter class
of spaces. 

One can develop homotopy and shape theories specifically designed
to work for at most $[L]$-dimensional spaces.
Absolute extensors for at most $[L]$-dimensional spaces
in a category of continuous maps are precisely $[L]$-soft mappings.
And compacta of trivial
$[L]$-shape are precisely $UV^{[L]}$-compacta~\cite{Ch}. One can define
(see \cite[Theorem 2.8]{Ch})
local $[L]$-contractibility in a standard way: a space $X$ is said to be
locally $[L]$-contractible (notation: $X \in LC^{[L]}$) if for any
neighbourhood $U$ of any point $x \in X$
there exists a smaller neighbourhood $V$ such that the inclusion
$V \hookrightarrow U$ is $[L]$-homotopic to a constant map. We present a full
proof (see Theorem \ref{dugun}) of a Dugunji-type theorem for such spaces.

We have several other applications of our results.
We characterize local $[L]$-softness
of a mapping in terms of local properties of the family
of its fibers (Theorem~\ref{thmcharLsoftness}).
This result was known for $n$-soft mappings~\cite{Dr1}.
Using idea from~\cite{B} on extension of $UV^n$-valued mappings,
we prove Theorem~\ref{thmextUVL} on extension of $UV^{[L]}$-valued
mappings. Also, we prove the following Theorem~\ref{thmfac} on factorization:
if the superposition $f\circ g$ of mappings of Polish spaces
is $[L]$-soft and $g$ is $UV^{[L]}$-map, then $f$ is $[L]$-soft.
For $n$-soft maps factorization theorem is proved in~\cite{B1}.

Another application is a version of Hurewicz
theorem for extension dimension. There are several approaches
to such a generalization of Hurewicz
theorem~\cite{DRS},\cite{ChV},\cite{Le},\cite{LL}.

\begin{A}
Let $f\colon X\to Y$ be a mapping of metric compacta where $\dim Y<\infty$.
Suppose that $\ed Y\leq [M]$ for some finite $CW$-complex~$M$.
If for some locally finite countable $CW$-complex $L$ we have
$\ed (f^{-1}(y)\times Z)\le [L]$ for every point $y\in Y$ and any
Polish space $Z$ with $\ed Z\le [M]$, then $\ed X\le [L]$.
\end{A}

The classical Hurewicz theorem for a mapping
$f \colon X \to Y$  of metric compacta with $\dim Y \leq m$ and
$\dim f = \sup\{ f^{-1}(y) \colon y \in Y\} \leq k$ follows from our result by letting $M = S^{m}$ and $L = S^{k+m}$. Indeed, note that 
$\dim\left( f^{-1}(y) \times Z \right) \leq k + m$ for any point $y \in Y$ and any
Polish space $Z$ with $\dim Z \leq m$. By our result, $\ed X \leq S^{k+m}$,
which means that $\dim X \leq k + m$ as required.

Section~\ref{S:Approx} of this paper is devoted to the
approximation theorem.
The {\it graph} of a multivalued mapping $F\colon X\to Y$
is the subset $\Gamma_F=\{(x,y)\in X\times Y\colon y\in F(x)\}$
of the product $X\times Y$.
We say that a multivalued mapping $F$ admits {\it approximations}
if every neighbourhood of the graph of $F$ contain the graph
of a singlevalued continuous mapping.

Usually one constructs
approximation as a composition of canonical mapping
into nerve of some covering and a mapping of this nerve,
defining the mapping of the nerve by induction on
dimension of its skeleta.
If the mapping is $UV^n$-valued and the domain space $X$ has
Lebesgue dimension $n$, then every point-image has trivial
shape relative to $X$ and relative to a nerve of some covering
of $X$, which allows one to construct a mapping from the nerve.
If extension dimension $\ed X=[L]$ does not coincide with
Lebesgue dimension of $X$, then $UV^{[L]}$-compactum
does not have trivial shape relative to a nerve of fine covering of $X$,
and one can not construct a mapping from the nerve.

Therefore, we have to define the approximation directly.
For some fine covering $\Sigma$ of $X$ we consider the sets
$\Sigma^{(k)}=\{x\in X\mid\ord_\Sigma x\le k+1\}$
and construct an approximation extending it successively
from $\Sigma^{(k)}$ to $\Sigma^{(k+1)}$.
Here $\Sigma^{(k)}$ plays a role of "$k$-dimensional skeleton"
of the cover $\Sigma$.
For elements $s_0,s_1,\dots,s_n\in\Sigma$ with non-empty
intersection $\cap_{i=0}^n s_i$ we consider the set
$\bigcup_{i=0}^n s_i\setminus\bigcup_{i\ne 0,1,\dots,n} s_i$
as a closed "simplex" with vertices $s_0,\dots,s_n$.
Also, we understand the set $\cap_{i=0}^n s_i$ as an interior of this simplex.
These notions of "skeleton" and "simplex" of a covering
allows us to proceed the proof in a usual way --- by induction
on "dimension" of "skeleta".
Note that our proof gives better result even for $UV^n$-valued mappings:
part (2) of Theorem~\ref{thmapproximation} was known only for metrizable
space $X$~\cite{Kr1}.

Sections~\ref{S:equi}--\ref{S:paracom} are devoted to selection problem.
The notion of filtration appeared to be very useful in continuous
selection theory (see~\cite{SB},~\cite{B2}) and we state our selection
theorem in terms of filtrations of multivalued mappings.

\begin{defin}
An increasing\footnote{We consider only increasing filtrations
indexed by a segment of the natural series starting from zero.}
finite sequence of subspaces
$$ Z_0\subset Z_1\subset\dots\subset Z_n\subset Z $$
is called a {\it filtration} of space $Z$ of length $n$. 
A sequence of multivalued mappings $\{F_k\colon X\to Y\}_{k=0}^n$
is called a {\it filtration of multivalued mapping} $F\colon X\to Y$
if $\{F_k(x)\}_{k=0}^n$ is a filtration of $F(x)$ for any $x\in X$.
\end{defin}

To construct a local selection we need our filtration of multivalued maps
to be complete and lower $[L]$-continuous.
The notion of completeness for multivalued mapping is introduced by
E.~Michael~\cite{Mi89}.

\begin{defin}
A multivalued mapping $G\colon X\to Y$ is called
{\it complete} if all sets $\{x\}\times G(x)$ are closed
with respect to some $G_\delta$-set $S\subset X\times Y$
containing the graph of this mapping.

We say that a filtration of multivalued mappings $G_i\colon X\to Y$
is {\it complete} if every mapping $G_i$ is complete.
\end{defin}

In section~\ref{S:equi} we introduce a notion of local property
of multivalued mapping. To have a local property, multivalued
mapping should have all fibers satisfying this local property,
and, moreover, the fibers should satisfy this property uniformly.
An important example of local property is local $[L]$-connectedness.

\begin{defin}\label{D:lconnected}
Let $L$ be a $CW$-complex. A pair of spaces $V\subset U$ is said
to be {\it $[L]$-connected}
if for every paracompact space $X$ of extension dimension $\ed X\le [L]$
and for every closed subspace $A\subset X$ any mapping of $A$
into $V$ can be extended to a mapping of $X$ into $U$.
\end{defin}

We call a multivalued mapping lower $[L]$-continuous
if it is locally $[L]$-connected:

\begin{defin}
A multivalued mapping $F\colon X\to Y$ is called
{\it $[L]$-continuous at a point} $(x,y)\in\Gamma_F$ of its graph if
for any neighbourhood $Oy$ of the point $y\in Y$, there are
a neighbourhood $O'y$ of the point $y$ and a neighbourhood $Ox$ of
the point $x\in X$ such that for all $x'\in Ox$, the pair
$F(x')\cap O'y\subset F(x')\cap Oy$ is $[L]$-connected.

A mapping which is $[L]$-continuous at all points of its graph is called
{\it lower $[L]$-continuous}.
We say that a filtration of multivalued mappings
is {\it lower $[L]$-continuous} if every mapping of this
filtration is lower $[L]$-continuous.
\end{defin}

To construct a global selection we need our filtration of multivalued maps
to be fiberwise $[L]$-connected.

\begin{defin}
A filtration of multivalued mappings $\{G_i\colon X\to Y\}_{i=0}^n$
is said to be {\it fiberwise $[L]$-connected} if for any point $x\in X$
and any $i<n$ the pair $G_i(x)\subset G_{i+1}(x)$ is $[L]$-connected.
\end{defin}

Now we can state our selection theorem.

\begin{B}
Let $L$ be a finite  $CW$-complex
such that $[L]\le [S^n]$ for some $n$.
Let $X$ be a paracompact space of extension dimension $\ed X\le [L]$.
Suppose that multivalued mapping $F\colon X\to Y$
into a complete metric space $Y$ admits a lower $[L]$-continuous, complete,
and fiberwise $[L]$-connected $n$-filtration
$F_0\subset F_1\subset\dots\subset F_n\subset F$.
If $f\colon A\to Y$ is a continuous singlevalued selection of $F_0$
over a closed subspace $A\subset X$, then there exists a continuous
singlevalued selection $\wt f\colon X\to Y$
of the mapping $F$ such that $\wt f|_A=f$.
\end{B}

Let us recall some definitions and introduce our notations.
We denote by $\Int A$ the interior of the set $A$.
For a cover $\omega$ of a space $X$ and for a subset $A\subseteq X$ let
$\St(A,\omega)$ denote the star of the set $A$ with respect to $\omega$.

For a subset $\mathcal U$ of the product $X\times Y$ we denote
by $\mathcal U(x)$ the subset $\pr_Y(\mathcal U\cap\{x\}\times Y)$ of $Y$,
where $x$ is a point of $X$.
For a multivalued mapping $F\colon X\to Y$ we denote by $F^\Gamma(x)$
the subset $\{x\}\times F(x)$ of $X\times Y$.
A multivalued mapping $F\colon X\to Y$
is said to be {\it upper semicontinuous} (shortly, u.s.c.) if its graph
is closed in the product $X\times Y$.
We say that multivalued mapping is {\it compact} if it is upper
semicontinuous and compact-valued.
A filtration consisting of compact
multivalued mappings is called {\it compact}.


A pair of subspaces $K\subset K'$ of a space $Z$
is called {\it $UV^{[L]}$-connected in $Z$}
if any neighbourhood $U$ of $K'$ contains a neighbourhood $V$ of $K$
such that the pair $V\subset U$ is $L$-connected.
A filtration $\{F_i\colon X\to Y\}_{i=0}^n$ of u.s.c. maps
is called {\it $UV^{[L]}$-connected $n$-filtration} if for any point $x\in X$
and any $i<n$ the pair $F_i(x)\subset F_{i+1}(x)$
is $UV^{[L]}$-connected in $Y$.
We say that multivalued mapping $F$ is {\it $n$-$UV^{[L]}$-filtered}
if it contains an $UV^{[L]}$-connected $n$-filtration.

A compact metric space $K$ is called {\it $UV^{[L]}$-compactum}
if the pair $K\subset K$ is $UV^{[L]}$-connected in any $ANR$-space.
Theorem~\ref{lemmainvar} shows that this property does not depend on
embedding of $K$ in Polish $ANE([L])$-space.
A multivalued mapping is called {\it $UV^{[L]}$-valued}
if it takes any point to $UV^{[L]}$-compactum.

A mapping $f\colon Y\to X$ is said to be {\it $[L]$-soft}
(resp. {\it locally $[L]$-soft})
if for any paracompact space $Z$ with $\ed Z\le [L]$, its closed subspace
$A\subset Z$ and any mappings $g\colon Z\to X$ and
$\wt g_A\colon A\to Y$ such that $f\circ\wt g_A=g|_A$
there exists a mapping $\wt g\colon Z\to Y$
(resp. $\wt g\colon OA\to Y$ of some neighbourhood of $A$)
such that $f\circ\wt g=g$ (resp. $f\circ\wt g=g|_{OA}$).
Finally let $AE([L])$ (resp. $ANE([L])$) denote the class of spaces with
$[L]$-soft (resp. locally $[L]$-soft) constant mappings.

\section{Singlevalued Approximation Theorem}\label{S:Approx}

We introduced in section~\ref{S:intro} the notions of "skeleton"
and "simplex" of a covering. For a covering $\Sigma$ of $X$
we denote by $\Sigma^{(k)}$ its $k$-dimensional skeleton
$\{x\in X\mid\ord_\Sigma x\le k+1\}$.
For elements $s_0,s_1,\dots,s_n\in\Sigma$ with non-empty
intersection $\cap_{i=0}^n s_i$ we define a
"closed $n$-dimensional simplex"
\[ [s_0,s_1,\dots,s_n]=\bigcup_{i=0}^n s_i\setminus
 \bigcup_{i\ne 0,1,\dots,n} s_i\]
and its "interior" $\langle s_0,s_1,\dots,s_n\rangle=\cap_{i=0}^n
s_i\cap\Sigma^{(n)}$.
It is easy to check that the $n$-skeleton consists of $n$-simplices
\[\Sigma^{(n)}=\bigcup\{[s_{i_0},s_{i_1},\dots,s_{i_n}]\mid
 \cap_{k=0}^n s_{i_k}\ne\emptyset\}\]
and that any "simplex" consists of its "boundary"
and its "interior"
\[ [s_0,s_1,\dots,s_n]=\bigcup_{m=0}^n [s_0,\dots,\wh s_m,\dots, s_n]
 \cup\langle s_0,s_1,\dots,s_n\rangle.\]
Clearly, $\Sigma^{(k)}$ is closed in $X$ and $\Sigma^{(n)}=X$
if the cover $\Sigma$ has order $n+1$.
The following property is important for our construction:
the "interiors" of distinct $k$-dimensional "simplices"
are mutually disjoint and
$$\Sigma^{(k)}=\bigcup\{\langle s_{i_0},s_{i_1},\dots,s_{i_n}\rangle\mid
 \cap_{k=0}^n s_{i_k}\ne\emptyset\}\cup\Sigma^{(k-1)}\eqno{(\dag)}$$

Suppose $Z$ is any space and $u$ is an open covering of $Z$.
We shall denote union of all elements of $u$ by $\cup u$.

Further we will consider triples of the form
$(X,\omega ,G)$, where $G$ is a multivalued mapping of $X$ to $Y$
and $\omega\in\cov X$.

\begin{defin} For a pair of spaces $X'\subset X$
a triple $(X',\omega', G')$ is said to be
$[L]$-connected refinement of a triple
$(X,\omega, G)$ if for any $W'\in\omega'$
there exists $W\in\omega$
with $\St (W',\omega')\subset W$ such that the pair $G' (\St (W',
\omega'))\subset G(W)$ is $[L]$-connected.

A sequence of triples $\{ (X_k ,\omega_k, G_k)\}_{k\le n}$ is said to be
$[L]$-connected if for each $k < n$ the triple $(X_k ,\omega_k, G_k)$ is
$[L]$-connected refinement of the triple $(X_{k+1},\omega_{k+1}, G_{k+1})$.
\end{defin}

\begin{lem} \label{sftsel}
Let $L$ be a $CW$-complex such that $[L]\le [S^n]$ for some $n$.
Let $X_0\subset\dots\subset X_{n+1}$ be a filtration of spaces
and $X$ be a paracompact subspace of a space $X_0$ such that $\ed X\le [L]$.

$(1)$ If $\{(X_k ,\omega_k, G_k)\}_{k\le n}$ is $[L]$-connected sequence
of triples, then there exists singlevalued continuous mapping
$f\colon X\to G_n(X_n)$ such that $f(x)\in G_n(\St(x,\omega_n))$
for each $x\in X$.

$(2)$ Suppose that $\{ (X_k ,\omega_k, G_k)\}_{k\le n+1}$
is $[L]$-connected sequence of triples.
Let $A$ be a closed subset of
$X$ and $g\colon A\to G_0 (X_0)$ be a singlevalued continuous mapping
such that $g(x)\in G_0 (\St (x,\omega_0))$ for each $x\in A$.
Then there exists singlevalued
continuous mapping $f\colon X\to G_{n+1} (X_{n+1})$
extending $g$ such that $f(x)\in G_{n+1}
(\St (x,\omega_{n+1}))$ for each $x\in X$.
\end{lem}

\begin{proof}
We shall prove the statement $(2)$. The proof of $(1)$ is similar.

Find an open locally finite covering $\Sigma$ of $X$
such that closures of elements of $\Sigma$ form strong star-refinement
of $\omega_0 |_X$ and order of $\Sigma$ is $\le n+1$.

Put $f_{-1} =g$. Let us construct
a sequence of mappings $\{ f_k\colon\Sgm{k}\bigcup A\to Y\}_{k=-1}^n$
such that $f_k$ extends $f_{k-1}$ and
$$f_k (x)\in G_{k+1} (\St (x,\omega_{k+1}))\mbox{ for each }
x\in\Sgm{k}\eqno{(*)}$$
Then we can let $f=f_n$ since $\Sigma^{(n)} =X$.

Suppose $f_k$ has been already constructed.
Since $(\dag)$ holds, it suffices to define $f_{k+1}$ on the "interior"
$\langle\sigma\rangle$ of each "simplex"
$[\sigma ]=[ s_0,s_1,\dots, s_{k+1} ]$.
Since $\Sigma$ is locally finite and the "interiors" of "closed
$k$-dimensional simplices" are mutually disjoint we can consider
each simplex independently.

Since $\omega_0$ is a star refinement of $\omega_{k+1}$,
there exists $V_{\sigma}\in\omega_{k+1}$ such
that $[\sigma ]\subset V_{\sigma}$
Since the triple $(X_{k+1} ,\omega_{k+1}, G_{k+1})$ is $[L]$-connected
refinement of the triple $(X_{k+2} ,\omega_{k+2}, G_{k+2})$, there exists
$U_{\sigma}\in\omega_{k+2}$ such that the pair
$G_{k+1}(\St(V_{\sigma},\omega_{k+1}))\subset G_{k+2}(U_{\sigma})$
is $[L]$-connected.

Let $[\sigma ]' = [\sigma ]\bigcap (A\bigcup\Sigma^{(k)})$.
For any $x\in [\sigma]'$ we have $x\in V_{\sigma}$
and the property $(*)$ implies
$f_k (x)\in G_{k+1} (\St (x,\omega_{k+1}))
\subset G_{k+1} (\St (V_{\sigma},\omega_{k+1}))$.
Hence $f_k ([\sigma ]')\subset G_{k+1} (\St
(V_{\sigma},\omega_{k+1}))$ and therefore $f_k$ can be extended over
$[\sigma ]$ to a map $\overline{f_k}\colon [\sigma ]\to
G_{k+2} (U_{\sigma})$.
We let $f_{k+1} |_{\langle\sigma\rangle} =
\overline{f_k}|_{\langle\sigma\rangle}$.

Let us check property $(*)$.
Since $\omega_{k+1}$ refines $\omega_{k+2}$, for all $x\in\Sigma
^{(k)}$ we have $f_{k+1} (x)=f_k (x)\in G_{k+1} (\St (x,\omega_{k+1}))
\subset G_{k+2} (\St (x,\omega_{k+2}))$. By $(\dag)$, any point
$x\in\Sigma^{(k+1)}\backslash\Sigma^{(k)}$ is contained in some
"interior" $\langle\sigma\rangle$. Since $\langle\sigma\rangle
\subset U_{\sigma}\in\omega_{k+2}$, we have
$f_{k+1}(x)\in G_{k+2} (U_{\sigma})\subset
G_{k+2} (\St (x,\omega_{k+2}))$.
\end{proof}

\begin{defin}
For a multivalued mapping $F\colon X\to Y$
an open neighbourhood $U\subset X\times Y$ of a fiber $F^{\Gamma(x)}$
is said to be {\it $F$-stable with respect to $x\in X$}
if there exists an open neighbourhood $O_x$ of
the point $x$ and an open subset $V_x\subset Y$
such that $\Gamma_{F|_{Ox}}\subset O_x\times V_x\subset U$.

The neighbourhood $U$ of the graph is said to be {\it $F$-stable}
if it is $F$-stable with respect to every point in $X$.
\end{defin}

\begin{defin}
A multivalued mapping $G\colon X\to Y$ is said to be
{\it a stable singular neighbourhood of $F$}
if for each $x\in X$ there exist open neighbourhoods $O_x$
of $x$ in $X$ and $V_x$ of $F(x)$ in $Y$ such that
$V_x\subset\bigcap\{ G(x')\mid x'\in Ox\}$.
\end{defin}

\begin{lem} \label{chain}
Let $X$ be a paracompact space and $L$ be a $CW$-complex.
Suppose that $\{ F_k\}_{k\le n}$ is a $UV^{[L]}$-connected
$n$-filtration consisting of multivalued mappings from $X$ to $Y$. Let
$\omega_n$ be a covering of $X$  and $G_n$ be a singular
stable neighbourhood of $F_n$. Then for each $k<n$ there exists an
open covering $\omega_k$ of $X$ and a
stable singular neighbourhood $G_k$ of
mapping $F_k$ such that the sequence $\{
(X,\omega_k, G_k)\}_{k\le n}$ is $[L]$-connected.
\end{lem}

\begin{proof}
We shall construct $\omega_k$ and $G_k$ by reverse
induction on $k$ starting from $k = n-1$. Since all inductive steps are
similar we shall show the constructions only for $k = n-1$.

Since $G_n$ is stable, for each $x\in X$ there exist
open neighbourhoods $O_x'$ of $x$ in $X$ and
$V_x'$ of $F_n (x)$ in $Y$ such that
$V_x'\subset\bigcap\{ G_n (x')\mid x'\in O_x'\}$.
Since $\{ F_k\}$ is $UV^{[L]}$-filtration there exist open neighbourhoods
$O_x\subset O_x'$ of $x$ and $V_x$ of $F_{n-1} (x)$ such that
$F_{n-1}(O_x)\subset V_x$ and the pair $V_x\subset V_x'$ is $[L]$-connected.
We may assume that the covering $\{ O_x\}_{x\in X}$ refines $\omega_n$.

Let $u\in\cov X$ be a locally finite strong star-refinement of
$\{ O_x\}_{x\in X}$.
For each $U\in u$ find $x(U)$ such that $\St (U,u)\subset O_{x(U)}$.
We shall also use notations $V_U =V_{x(U)}$ and $O_U =O_{x(U)}$.

For each $x\in X$ we put $G_{n-1} (x)=\bigcap\{ V_U\mid x\in U\}$.
Let $\omega_{n-1}$ be a strong star-refinement of $u$.

Let us check that $G_{n-1}$ is a stable singular neighbourhood of
$F_{n-1}$. Consider any $x\in X$. Find open neighbourhood $Ox$ of $x$
which intersects only finitely many elements of $u$. We may assume that
$Ox\subset Ux$ for some $Ux\in u$.
Put $Vx=\bigcap\{ V_U\mid U\cap Ox\ne\varnothing\}$.
Since $x\in Ox\subset Ux$ it follows by the choice of $O_U$ that for all
$U$ such that $U\cap Ox\ne\varnothing$
we have $x\in O_U$. Hence, using the fact $F_{n-1} (O_U)\subset V_U$
we obtain $F_{n-1} (x)\subset Vx$.
Finally, we have $\bigcap\{ G(x')\mid x'\in Ox\}=\bigcap\{\bigcap\{ V_U\mid x'\in U\}
\mid x'\in Ox\} =Vx$ by the definition of $Vx$.

Let us show that $(X,\omega_{n-1} ,G_{n-1})$ is $[L]$--connected
refinement of the triple $(X,\omega_{n} ,G_{n})$.
Consider any $W'\in\omega_{n-1}$. Find $U'\in u$ such that
$\St (W',\omega_{n-1})\subset U'$.
There exists $W\in\omega_n$ with $O_{U'}\subset W$.
Take $x\in\St (W',\omega_{n-1})$. Then $G_{n-1} (x) =\bigcap\{
V_U\mid x\in U\}\subset V_{U'}$ and the pair $V_{U'}\subset V'_{x(U')}$
is $[L]$--connected. Finally, observe that by the choice of $\{ O_x'\}$ and
$\{ V_x'\}$ we have $V'_{x(U')}\subset\bigcap\{ G_n (x')\mid x'\in
O_{U'}\}\subset G_n (W)$.
\end{proof}

\begin{thm} \label{thmapproximation}
Let $L$ be a $CW$-complex such that $[L]\le [S^n]$ for some $n$.
Let $X$ be a paracompact space of extension dimension $\ed X\le [L]$.

$(1)$ If $F\colon X\to Y$ is a multivalued mapping which admits
$UV^{[L]}$-connected $n$-filtration, then any $F$-stable neighbourhood of
the graph $\Gamma_{F}$ contains a graph of a singlevalued continuous
mapping of $X$ to $Y$.

$(2)$Let $A\subset X$ be a closed subspace.
If $F$ admits $UV^{[L]}$-connected $(n+1)$-filtration
$F_0\subseteq F_1\subseteq\dots\subseteq F_{n+1}$, then for any
$F$-stable neighbourhood $U$ of the graph $\Gamma_{F}$ there exists
$F_0$-stable neighbourhood $V$ of the graph $\Gamma_{F_0|_A}$
such that every singlevalued continuous mapping $g\colon A\to Y$
with $\Gamma_{g}\subset V$ can be extended to a singlevalued
continuous mapping $f\colon X\to Y$ with $\Gamma_{f}\subset U$.
\end{thm}

\begin{proof}
We shall prove statement $(2)$. The proof of $(1)$ is similar.
Let $U$ be an arbitrary stable neighbourhood of the graph of $F$.
Since $U$ is stable, for each $x\in X$ there
exist open neighbourhoods $O_x$ of $x$
and $V_x$ of $F(x)$ such that
$\Gamma_{F|_{O_x}}\subset O_x\times V_x\subset U$.
Let $\omega_{n+1}$ be a strong star refinement of
$\{ O_x\}_{x\in X}$.

For each $x\in X$ we let
$G_{n+1} (x)=\bigcap\{ U(x')\mid x'\in\St (x,\omega_{n+1})\} $.
Let us check that $G_{n+1}$ is a stable singular neighbourhood of $F_n$.
Fix $x\in X$ and consider $W\in\omega_{n+1}$ which contains $x$. Then
$$\bigcap\{ G_{n+1} (x')\mid x'\in W\} =
\bigcap\{\bigcap\{ U(x'')\mid x''\in\St (x' ,\omega_{n+1})\}
\mid x'\in W\}$$
$$\supset\bigcap\{\bigcap\{ U(x'')\mid x''\in\St (W,\omega_{n+1})\}
\mid x'\in W\}\supset V_z\supset F(x)$$
where $z\in X$ is chosen so that $\St(W,\omega_{n+1})\subset O_z$.

Using Lemma~\ref{chain}, construct an $[L]$-connected
sequence $\{ (X,\omega_k, G_k)\}_{k\le n+1}$.
Observe that since $G_0$ is stable singular neighbourhood of $F_0$,
the graph $\Gamma_{G_0}$ contains an open stable neighbourhood $V$
of $\Gamma_{F_0}$.

Suppose that $g\colon A\to Y$ is a singlevalued continuous mapping
such that graph of $g$ is contained in $V$.
Then $g(x)\in G_0 (x)$ for all $x\in A$.
Hence we can apply Lemma~\ref{sftsel} and obtain singlevalued continuous
mapping $f\colon X\to Y$ extending $g$ such
that $f(x)\in G_{n+1} (\St (x,\omega_{n+1}))$ for each $x\in X$.
This fact and the definition of $G_{n+1}$ imply that graph of $f$ is
contained in $U$.
\end{proof}

\begin{lem} \label{chainmetr}
Let $X$ be a subspace of a metric space $M$ and
$\mathcal U_n$ be an open neighbourhood of $X$ in $M$.
For a $CW$-complex $L$ suppose that $\{ F_k\colon X\to Y\}_{k\le n}$
is a $UV^{[L]}$-connected $n$-filtration.
Let $\omega_n$ be a covering of $\mathcal U_n$ and
$G_n\colon\mathcal U_n\to Y$ be a stable singular neighbourhood of $F_n$.
Then there exists $[L]$-connected sequence
$\{(\mathcal U_k ,\omega_k, G_k)\}_{k\le n}$ such that
$\mathcal U_k$ is an open neighbourhood of $X$ in $M$ and
$G_k$ is a stable singular neighbourhood of $F_k$.
\end{lem}

\begin{proof}
We shall construct $\mathcal U_k$, $\omega_k$ and $G_k$ by reverse
induction on $k$ starting from $k = n-1$. Since all inductive steps are
similar we shall show the constructions only for $k = n-1$.

Since $G_n$ is stable, for each $x\in X$ there exist
open neighbourhoods $O_x'$ of $x$ in $U_n$ and
$V_x'$ of $F_n (x)$ in $Y$ such that
$V_x'\subset\bigcap\{ G_n (x')\mid x'\in O_x'\}$.
Since $\{ F_k\}$ is $UV^{[L]}$-filtration
there exist open in $M$
neighbourhood $O_x\subset O_x'$ of $x$ and open neighbourhood
$V_x$ of $F_{n-1} (x)$
such that $F_{n-1} (O_x)\subset V_x$ and
the pair $V_x\subset V_x'$
is $[L]$-connected.
We may assume that the collection $\{ O_x\}_{x\in X}$
refines $\omega_n$. Put $\mathcal U_{n-1} =\bigcup\{ O_x\mid x\in X\}$.

Let $u$ be a locally finite covering of $\mathcal U_{n-1}$
which is a strong star-refinement of $\{ O_x\}_{x\in X}$.
For each $U\in u$ find $x(U)$ such that $\St (U,u)\subset O_{x(U)}$.
For any $x\in\mathcal U_{n-1}$ we put $G_{n-1} (x)=\bigcap\{ V_{x(U)}\mid x\in U\}$.
Let $\omega_{n-1}$ be a strong star-refinement of $u$.
Then similarly to the proof of Lemma~\ref{chain}
we obtain that $G_{n-1}$ is a stable singular neighbourhood of
$F_{n-1}$ and the triple $(\mathcal U_{n-1},\omega_{n-1} ,G_{n-1})$
is $[L]$-connected refinement of the triple $(\mathcal U_n ,\omega_n ,G_n)$
\end{proof}

\begin{defin}
A singlevalued continuous
surjective mapping $f\colon Y\to X$ of metric spaces is said to
be approximately $[L]$-invertible if for any embedding of $f$ into the
projection $p\colon M\times N\to M$ of metric spaces
where $M\in ANE([L])$
the following condition is
satisfied:

for any neighbourhood $W$ of $Y$ in $M\times N$ there exists open
neighbourhood $U$ of $X$ in $M$ such that for any mapping
$g\colon Z\to U$ of
paracompact space $Z$ with $\ed (Z)\le [L]$ there exists a lifting $g'
\colon Z\to W$ of $g$ such that $pg' = g$.
\end{defin}

\begin{thm}\label{thmapprinvert}
Let $L$ be a $CW$-complex such that $[L]\le [S^n]$ for some $n$.
Suppose that for a continuous singlevalued surjective mapping of
metric spaces $f$ the multivalued mapping $F=f^{-1}$
admits a compact $UV^{[L]}$-connected $n$-filtration.
Then $f$ is approximately $[L]$-invertible.
\end{thm}

\begin{proof}
Consider an embedding of $f$ into the projection $p\colon M\times N\to M$
of metric spaces where $M\in ANE([L])$ and fix an arbitrary neighbourhood
$W$ of $Y$ in $M\times N$.
Let $\{ F_i\}_{i=0}^n$ be a compact
$UV^{[L]}$--connected $n$-filtration of $F=f^{-1}$.
Then the mapping $F' = pr_N\circ F$ admits a compact
$UV^{[L]}$--connected $n$-filtration $\{ F_i =pr_N\circ F_i\}_{i=0}^n$.

Since the mapping $F_n'$ is compact,
$W$ is a stable neighbourhood of the
graph $\Gamma_{F_n'}\subset M\times N$.

For each $x\in X$ find open neighbourhood $O_x$ of $x$ in $M$ and open
subset $V_x$ of $N$ such that
$\Gamma_{F_n' |_{O_x}}\subset O_x\times V_x\subset W$.
Let $\mathcal U_n =\bigcup\{ O_x\mid x\in X\}$ and $\omega_n\in\cov U_n$
be a strong star refinement of $\{ O_x\}_{x\in X}$. We can define a stable
singular neighbourhood $G_n$ of $F_n'$ letting, as before,
$G_n (x)=\bigcap\{W(x')\mid x'\in\St (x,\omega_n)\}$
for all $x\in\mathcal U_n$.
By Lemma~\ref{chainmetr} we can find $[L]$-connected sequence of triples
$\{ (\mathcal U_k ,\omega_k ,G_k)\}_{k\le n}$
where $G_k\colon\mathcal U_k\to N$ is a
stable singular neighbourhood of $F_k'$.

Put $U=\mathcal U_0$ and show that the pair $(W,U)$ satisfies lifting property.
Consider an arbitrary mapping $g\colon Z\to U$ where $Z$ is
a paracompact space with $\ed Z\le [L]$.
We may assume that $g$ is embedded into a projection $p'\colon M\times E\to M$
for some Tychonov space $E$ such that $Z\subset M\times E$.
For each $k = 0,1,\dots,n$ we let $\mathcal U_k' = (p')^{-1}\mathcal U_k$
and define open in $M\times E$ covering $\omega_k' = (p')^{-1}\omega_k$
of $\mathcal U_k'$ and multivalued mapping $G_k'\colon U_k'\to N$ letting
$G'_k (x) = G_k (p' (x))$ for all $x\in\mathcal U_k'$.
It is easily seen that the sequence
$\{(\mathcal U_k',\omega'_k ,G_k')\}_{k\le n}$ is also $[L]$-connected.
Hence we can apply Lemma~\ref{sftsel} to obtain a map $h\colon Z\to N$
such that $h(z)\in G_n'(\St (z,\omega'_n))$ for all $z\in Z$.

Now we can define lifting map $g'$ on $Z$ letting $g' (z) = (g(z),
h(z))$.
Clearly $pg' = g$. It is easel seen from the construction and
definition of $G_n$ that $g'$ maps $Z$ into $W$.
\end{proof}

\section{Local properties of multivalued mappings}\label{S:equi}

We follow definitions and notations from~\cite{DM}.

\begin{defin}
An ordering $\alpha$ of the subsets of a space $Y$ is {\it proper} provided:
\begin{itemize}
\item[(a)] If $W\alpha V$, then $W\subset V$;
\item[(b)] If $W\subset V$, and $V\alpha R$, then $W\alpha R$;
\item[(c)] If $W\alpha V$, and $V\subset R$, then $W\alpha R$.
\end{itemize}
\end{defin}

Further we will not mention the space on which the proper ordering
is defined.

\begin{defin}
Let $\alpha$ be a proper ordering.
\begin{itemize}
\item[(a)]
A metric space $Y$ is {\it locally of type $\alpha$} if, whenever
$y\in Y$ and $V$ is a neighbourhood of $y$, then there a neighbourhood
$W$ of $y$ such that $W\alpha V$.
\item[(b)]
A multivalued mapping $F\colon X\to Y$ of topological space $X$
into metric space $Y$ is {\it lower $\alpha$-continuous}
if for any points $x\in X$ and $y\in F(x)$ and for any neighbourhood $V$
of $y$ in $Y$ there exist neighbourhoods $W$ of $y$ in $Y$
and $U$ of $x$ in $X$ such that $(W\cap F(x'))\alpha(V\cap F(x'))$
provided $x'\in U$.
\end{itemize}
\end{defin}

For example, if $W\alpha V$ means that $W$ is contractible in
$V$, then locally of type $\alpha$ means locally contractible.
Another topological property which arise in this manner
is $LC^n$ (where $W\alpha V$ means that every continuous mapping
of the $n$-sphere into $W$ is homotopic to a constant mapping in $V$).
For the special case $n=-1$ the property $W\alpha V$ means that $V$ is
non-empty, and lower $\alpha$-continuity is lower semicontinuity.

If $W\alpha V$ means that the pair $W\subset V$ is $[L]$-connected,
then locally of type $\alpha$ means local absolute extensor
in dimension $[L]$.
And we call lower $\alpha$-continuity of multivalued mapping
as lower $[L]$-continuity.

\begin{lem} \label{lemmaequipoint}
Let $F\colon X\to Y$ be lower $\alpha$-continuous multivalued mapping
of topological space $X$ to metric space $Y$.
Consider a point $y\in F(x)$.
Then for any $\ep>0$ there exist $\delta>0$ and neighbourhoods $O_y$ of
the point $y$ in $Y$ and $O_x$ of the point $x$ in $X$ such that for
any points $x'\in O_x$ and $y'\in F(x')\cap O_y$ we have
$(O(y',\delta)\cap F(x'))\alpha (O(y',\ep)\cap F(x'))$.
\end{lem}

\begin{proof}
Since the mapping $F$ is lower $\alpha$-continuous, there are positive
$\delta<\ep/4$ and a neighbourhood $O_x$ of the point $x$ such that
$(O(y,2\delta)\cap F(x'))\alpha (O(y,\ep/2)\cap F(x'))$
for every point $x'\in O_x$.
Put $O_y=O(y,\delta)$. Then for every $x'\in O_x$ and every
$y'\in F(x')\cap O_y$ we have inclusions
$O(y',\delta)\subset O(y,2\delta)$ and $O(y,\ep/2)\subset O(y',\ep)$.
Therefore, $(O(y',\delta)\cap F(x'))\alpha (O(y',\ep)\cap F(x'))$.
\end{proof}

\begin{lem} \label{lemmaequicompactset}
Let $F\colon X\to Y$ be lower $\alpha$-continuous multivalued mapping
of topological space $X$ to metric space $Y$.
Consider a compact subset $K$ of the fiber $F(x)$.
Then for any $\ep>0$ there exist $\delta>0$ and neighbourhoods $OK$ of
compactum $K$ in $Y$ and $O_x$ of the point $x$ in $X$ such that for
any points $x'\in O_x$ and $y'\in F(x')\cap OK$ we have
$(O(y',\delta)\cap F(x'))\alpha (O(y',\ep)\cap F(x'))$.
\end{lem}

\begin{proof}
For every point $y\in K$ take a number $\delta_y>0$ and neighbourhoods
$Oy$ of the point $y$ and $O_yx$ of the point $x$ by Lemma~\ref{lemmaequipoint}.
Choose a finite subcovering $\{Oy_i\}_{i=1}^m$ of the cover
$\{Oy\}_{y\in K}$ of compactum $K$ and consider the corresponding
numbers $\delta_1,\dots,\delta_m$ and neighbourhoods
$O_1x,\dots,O_mx$ of the point $x$.
Clearly, we can put
$$ OK=\bigcup_{i=1}^m Oy_i, \qquad \delta=\min_{1\le i\le m}\delta_i,
\qquad Ox=\bigcap_{i=1}^m O_ix. $$
The lemma is proved.
\end{proof}

\begin{lem} \label{lemmaequicompactmap}
Suppose that lower $\alpha$-continuous multivalued mapping $F\colon X\to Y$
of paracompact space $X$ to metric space $Y$
contains a compact submapping $H\colon X\to Y$.
Then for any continuous positive function $\ep\colon X\to\R$
there exist a continuous positive function $\delta\colon X\to\R$
and a neighbourhood $U$ of the graph $\Gamma_H$ such that
for any points $x\in X$ and $y\in F(x)\cap U(x)$ we have
$(O(y,\delta(x))\cap F(x))\alpha (O(y,\ep(x))\cap F(x))$.
\end{lem}

\begin{proof}
Using Lemma~\ref{lemmaequicompactset}, we can find for
every point $x\in X$ a number
$\sigma(x)$ and open neighbourhoods $Ox$ of the point $x$
and $OH(x)$ of the compactum $H(x)$ such that
$(O(y',\sigma(x))\cap F(x'))\alpha (O(y',\ep(x)/2)\cap F(x'))$
for any points $x'\in Ox$ and $y'\in F(x')\cap OH(x)$.
Moreover, we may take a neighbourhood $Ox$ to be so small
that $H(Ox)$ is contained in $OH(x)$ and
$\sup_{x'\in Ox}\ep(x')<2\cdot\inf_{x'\in Ox}\ep(x')$.

Let us refine a locally finite cover
$\omega=\{W_\lambda\}_{\lambda\in\Lambda}$
into the cover $\{Ox\}_{x\in X}$ and for every $\lambda\in\Lambda$
take a point $x_\lambda$ such that $W_\lambda$ is contained in $Ox_\lambda$.
Let $\delta\colon X\to\R$ be a continuous positive function
such that for every point $x\in X$ we have
$\delta(x)\le\min\{\sigma(x_\lambda)\mid x\in W_\lambda\}$.
Put $U=\cup_{\lambda\in\Lambda} W_\lambda\times OH(x_\lambda)$.
Since $H(W_\lambda)$ is contained in $OH(x_\lambda)$ and
the sets $W_\lambda$ cover $X$, then $U$ is a neighbourhood of
the graph $\Gamma_H$.

Consider an arbitrary point $\{x\}\times\{y\}\in U\cap\Gamma_F$.
By the construction of $U$, there is a set $W_\lambda$ containing $x$
such that $\{x\}\times\{y\}\in W_\lambda\times OH(x_\lambda)$. Then
$(O(y,\sigma(x_\lambda))\cap F(x))\alpha (O(y,\ep(x_\lambda)/2)\cap F(x))$.
Therefore, since $\varepsilon(x)>\varepsilon(x_\lambda)/2$ and
$\delta(x)\le\sigma(x_\lambda)$, we have
$(O(y,\delta(x))\cap F(x))\alpha (O(y,\ep(x))\cap F(x))$.
\end{proof}

In what follows we are going to work with covers
of the product $X\times Y$ of paracompact space $X$ and metric space $Y$.
It will be convenient to work with "rectangular" covers.
And we consider covers of the form $\omega\times\ep$
where $\omega$ is a covering of $X$ and $\ep\colon X\to\R$
is a continuous positive function. Precisely, the covering
$\omega\times\ep$ consists of all products
$\{ W\times O(y,\varepsilon(x))|\; x\in W\in\omega,\ x\in X\}$.

\begin{rem}
A real-valued function $\varepsilon\colon X\to \R$ is called {\it locally
positive} if for any point $x$, there exists a neighbourhood
on which the infimum of the function is positive. For any locally
positive function $\varepsilon(x)$ on a paracompact space, there
exists a positive continuous function which is less than this function.
Indeed, consider a partition of the unity
$\{\varphi_\alpha(x)\}$ subordinated to a locally finite covering
$\{W_\alpha\}$ of this paracompact space where
the function $\varepsilon(x)$ is greater than
some positive number $c_\alpha$ on each element $W_\alpha$
of this covering. Then the function $\sum_\alpha
c_\alpha\cdot\varphi_\alpha(x)$ is the desired continuous function.
\end{rem}

The following lemma shows that if we have a graph
$\Gamma_H\subset X\times Y$ of a compact multivalued mapping
$H\colon X\to Y$
of paracompact space $X$ to metric space $Y$,
then we may consider only "rectangular" covers
of this graph of the form $\omega\times\ep$.

\begin{lem} \label{lemmarectcover}
For any open cover $\gamma$ of the graph $\Gamma_H\subset X\times Y$
of a compact multivalued mapping $H\colon X\to Y$
of paracompact space $X$ to metric space $Y$
there exist
an open cover $\omega\in\cov X$ and a continuous positive function
$\ep\colon X\to \R$ such that the cover $\omega\times\varepsilon$
of the graph $\Gamma_H$ refines $\gamma$.
\end{lem}

\begin{proof}
Consider a point $x\in X$.
For every point $\{x\}\times \{y\}\in \{x\}\times H(x)$
we fix its open neighbourhood $O_yx\times Oy$ refining $\gamma$.
Take a finite subcover $\{Oy_i\}_{i=1}^N$ of the cover
$\{Oy\}_{y\in H(x)}$ of the compactum $H(x)$ and let $2\lambda (x)$
be its Lebesgue number. We put
$$ Ox=(\bigcap_{i=1}^N O_{y_i}x)\cap \{x'\in X\mid H(x')\subset
  O(H(x),\lambda(x))\} $$
Then for any points $x'\in Ox$ and $y'\in H(x')$ the set
$Ox\times O(y',\lambda(x))$ refines $\gamma$.
Consider an open locally finite cover $\omega\in\cov X$
refining the cover $\{Ox\}_{x\in X}$.
For every $W\in \omega$ we fix an element $Ox_W$ of the
cover $\{Ox\}_{x\in X}$ such that $W\subset Ox_W$.
Since the cover $\omega$ is locally finite, the function
$\varepsilon'(x)=\min\limits_{x\in W\in\omega} \lambda(x_W)$
is locally positive. Let $\ep$ be any positive continuous function
which is less than $\ep'$. Then we define
$ \omega\times\varepsilon=
 \{W\times O(y,\varepsilon(x))\mid x\in W\in\omega, y\in H(x)\subset Y\}$.
\end{proof}

In what follows we shall construct for a given positive
continuous function $\delta\colon X\to \R$ an open covering
$\omega\in\cov X$ such that the function $\delta$ vary within
any element of the covering $\omega$ less than by half
(i.e. $\sup_{x\in W}\delta(x)<2\cdot\inf_{x\in W}\delta(x)$).
The following lemma shows the reason for such construction.

\begin{lem} \label{lemmavarybyhalf}
Suppose that a positive continuous function $\delta\colon X\to \R$ vary
within any element of the covering $\omega\in\cov X$ less than by half.
Then for any points $p_0=\{x_0\}\times\{y_0\}\in X\times Y$ and
$p=\{x\}\times\{y\}\in\St(p_0,\omega\times\delta)$
the star $\St(p_0,\omega\times\delta)$ is contained in the product
$\St(x_0,\omega)\times O(y,16\cdot\delta(x))$.
\end{lem}

\begin{proof}
For any point $x'\in \St(x_0,\omega)$ we have
$\delta(x')\le 2\cdot\delta(x_0)\le 4\cdot\delta(x)$.
Then the distance between points $y_0$ and $y$ is less than $8\cdot\delta(x)$.
Clearly, every element of the cover $\omega\times\delta$ containing
the point $p_0$ lies in the set
$\St(x_0,\omega)\times O(y_0,8\cdot\delta(x))$.
Therefore, the star $\St(p_0,\omega\times\delta)$
is contained in the product
$\St(x_0,\omega)\times O(y,\dist(y,y_0)+8\cdot\delta(x))$.
The lemma is proved.
\end{proof}

Let a lower semicontinuous mapping $\Phi\colon X\to Y$
contain a compact submapping $\Psi$.
Let us define the notion of starlike $\alpha$-refinement,
relative to a pair $(\Psi,\Phi)$, of coverings of the form
$(\omega\times\varepsilon)$, where $\omega\in \cov X$ and
$\varepsilon$ is a positive continuous function on $X$.

\begin{defin}
A covering $(\omega'\times\varepsilon')$ is called
{\it starlike $\alpha$-refined} into a covering
$\omega\times\varepsilon$ relative to a pair $(\Psi,\Phi)$ if
for any point $z\in \St(\Gamma_\Psi,\omega'\times\varepsilon')$
there exists an element $W\times O(y,\ep(x))$ of the cover
$\omega\times\varepsilon$ containing the star
$\St(z,\omega'\times\varepsilon')$ and such that
$$(\St(z,\omega'\times\varepsilon')(x')\cap\Phi(x'))\alpha
(O(y,\ep(x))\cap\Phi(x'))$$
for any point $x'\in \pr_X(\St(z,\omega'\times\varepsilon'))$.
\end{defin}

\begin{lem} \label{lemmastarlikerefinement}
Suppose that lower $\alpha$-continuous multivalued mapping $F\colon X\to Y$
of paracompact space $X$ to metric space $Y$
contains a compact submapping $H\colon X\to Y$.
Then for any continuous positive function $\ep\colon X\to \R$
and any open cover $\omega\in cov X$
there exist a continuous positive function $\delta\colon X\to \R$
and an open cover $\omega'\in cov X$ such that the cover
$\omega'\times\delta$ is starlike $\alpha$-refined
into a covering $\omega\times\varepsilon$ relative to a pair $(H,F)$.
\end{lem}

\begin{proof}
By Lemma~\ref{lemmaequicompactmap}
there exist a neighbourhood $\mathcal U$ of the graph $\Gamma_H$
and continuous positive function $\sigma\colon X\to \R$ such that
$16\sigma<\ep$ and for any points $x\in X$ and $y\in F(x)\cap\mathcal U(x)$
we have $(O(y,16\sigma(x))\cap F(x))\alpha(O(y,\ep(x))\cap F(x))$.
By Lemma~\ref{lemmarectcover} there is a covering
$\omega''\times\nu$ of the graph $\Gamma_H$ such that the star
$\St(\Gamma_H,\omega''\times\nu)$ is contained in $\mathcal U$.
Define a continuous positive function $\delta\colon X\to \R$
by the equality $\delta(x)={\frac{1}{16}}\min\{\sigma(x),\nu(x)\}$.
Consider a covering $\omega'\in cov X$ which is starlike refined
into $\omega$ and $\omega''$ and such that the function $\ep$ vary within
any element of the covering $\omega'$ less than by half.

Then for every point
$p_0=\{x_0\}\times\{y_0\}\in\St(\Gamma_H,\omega'\times\delta)$
the star $\St(p_0,\omega'\times\delta)$ is contained in $\mathcal U$.
Indeed, the star $\St(x_0,\omega')$ is contained in some element
$V$ of the cover $\omega''$. Take a point
$p=\{x\}\times\{y\}\in\Gamma_H\cap\St(p_0,\omega'\times\delta)$.
By the construction of the cover $\omega''\times\nu$
the set $V\times O(y,\nu(x))$ is contained in $\mathcal U$.
By Lemma~\ref{lemmavarybyhalf} the star $\St(p_0,\omega'\times\delta)$
is contained in $V\times O(y,16\delta(x))$.

Consider an arbitrary point $x'\in \St(x_0,\omega')$ and suppose
that the intersection of the set $\St(p_0,\omega'\times\delta)(x')$
with the fiber $F(x')$ is not empty and contains a point $y'$.
Then this intersection is contained in $O(y',16\delta(x'))$.
Since the point $\{x'\}\times\{y'\}$ lies in $\mathcal U$, then
$(O(y',16\sigma(x'))\cap F(x))\alpha (O(y',\ep(x'))\cap F(x'))$.
Fix an element $W$ of the cover $\omega$ containing the star
$\St(x_0,\omega')$. Clearly, the element $W\times O(y',\ep(x'))$
of the cover $\omega\times\varepsilon$ contains the star
$\St(p_0,\omega'\times \delta)$ (we apply Lemma~\ref{lemmavarybyhalf})
and the set $\{x'\}\times O(y',\ep(x'))$.
\end{proof}

The set
$$\st(A,\omega)=\bigcup\{U\in\omega|\; A\subset U\}$$
is the {\it small star} of a set $A$ relative to a covering $\omega$.
The proof of the following lemma is easy (actually, it is
Lemma of Continuity of Star Trace from~\cite{SB})

\begin{lem} \label{lemmastartrace}
Let $\omega$ be an open covering of a metric space $Y$, let
$F\colon X\to Y$ be a compact multivalued mapping, and let
$\Phi\colon X\to Y$ be complete lower $\alpha$-continuous mapping.
Then the multivalued mapping $G$ which assigns the set
$\Phi(x)\cap \st(F(x),\omega)$ to the point $x\in X$
is complete and lower $\alpha$-continuous.
\end{lem}

\begin{proof}
The multivalued mapping $G'$ which assigns the small star $\st(F(x),\omega)$
to a point $x\in X$ has the open graph in the space $X\times Y$.
Indeed, for a point $\{x\}\times \{y\}\in \Gamma_{G'}$
there is an element $W\in \omega$ containing the image $F(x)$.
Then by the upper semicontinuity of $F$, for some neighbourhood
$Ox\subset X$ of the point $x$, the image $F(Ox)$ is contained in $W$.
Then the set $Ox\times W$ is an open neighbourhood of the point
$\{x\}\times\{y\}$ in the graph $\Gamma_{G'}$.

Now the completeness and the lower $\alpha$-continuity of mapping $\Phi$
imply these properties for the mapping $G=G'\cap\Phi$ by the openness
of the graph $\Gamma_{G'}$.
\end{proof}

\section{$[L]$-soft mappings}\label{S:soft}

In this section we prove several important technical results
about $[L]$-soft mappings. In particular, these results
allows us to show that $UV^{[L]}$-property of compactum
does not depend on embedding of this compactum into $ANE([L])$-space.

\begin{thm}\label{dugun}
Let $L$ be a locally finite countable $CW$-complex such that
$[L]\le [S^n]$ for some $n$.
Then for a Polish space $Y$ property $Y\in LC^{[L]}$ implies $Y\in ANE([L])$.
\end{thm}
\begin{proof}
By Proposition \ref{polishane}, it suffices to check property
$Y\in ANE([L])$ for Polish spaces.
Since any Polish space $X$ with $\ed X\le [L]$ admits
closed embedding into Polish $AE([L])$-space of extension dimension
$\le [L]$~\cite{Ch97}, we may assume that $X \in AE([L])$.

Let $A$ be a closed subspace of $X$ and $f \colon A\to Y$ be
a continuous mapping.
There is an open covering $\omega$ of $X\backslash A$
with the following property: $(i)$ for any point $a\in A$ and any its
neighbourhood $O_a$ in $X$ there exists a neighbourhood $V_a$ of
$a$ in $X$ such that for all $W\in \omega$ if
$W\cap V_a \ne \varnothing$ then $U\subset O_a$~\cite[Theorem~3.1.4]{Bor}.
Since $\dim (X\backslash A) \le n$ there exists an open refinement
$u=\bigcup_{k=0}^n u_k$ of $\omega$ where $u_k$ is a countable
discrete system of open disjoint sets~\cite{Eng}.

For each $U^0_i \in u_0$ choose $a_i \in A$ such that
$\dist (a_i,U^0_i) \le \sup\{\dist (x,A) \mid x\in U^0_i\}$
and define a mapping $f_0$ on $W_0=\cup \{U^0_i\mid U^0_i \in u_0\}\cup A$
as follows: $f_0|_A = f|_A$ and $f_0 (U^0_i) = f(a_i)$.
It is easily seen that $f_0$ is continuous.

By induction on $k=1,\dots,n$ we shall
find neighbourhoods $W_k$ of $A$ in
$\bigcup_{j=0}^k \{ U^j_i \mid U^j_i \in u_j\} \bigcup A$
and using $f_{k-1}$ we shall extend $f$ to $f_k \colon W_k \to Y$.
Since $u$ covers $X\backslash A$ the mapping $f_n$
extends $f$ to the neighbourhood $W_n$ of $A$ in $X$.

Suppose that $f_{k-1}$ has been already constructed.
Since $Y \in LC^{[L]}$, for each $a\in A$ there exists a neighbourhood
$O_a$ of $a$ in $X$ such that $f_{k-1} |_{O_a}$ is $[L]$-homotopic to a
constant map in $Y$.
Applying to $O_a$ property $(i)$ of $u$ find
neighbourhood $V_a \subset O_a$.
Put ${\cal V}_k =\bigcup \{ V_a \mid a\in A\}$ and
$W_k =\cup \{ U^k_i \mid U^k_i \subset {\cal V}_k \} \cup W_{k-1}$.
Observe that for all $U^k_i \in u_k$ we have:
$(ii)$ $f_{k-1} |_{ U^k_i \cap W_{k-1}}$ is $[L]$-homotopic to a
constant map in $Y$ provided $U^k_i \subset {\cal V}_k$.

We shall define $f_k$ as an extension of $f_{k-1}$ from the set
$W_{k-1} \backslash (\cup \{ U^k_i \mid U^k_i \subset {\cal V}_k \})$.
Since the system $u_k$ is disjoint, we can
define $f_k$ independently on every $U^k_i \subset {\cal V}_k$.
Consider an arbitrary $U^k_i \in u_k$ such that $U^k_i \subset {\cal V}_k$.
If $W_{k-1}\backslash U^k_i$ is open in $X$, choose a point $a_i \in A$
such that $\dist (a_i,U^k_i)\le \sup \{ \dist (x,A)\mid x\in U^k_i \}$
and define $f_k (U^k_i)=f(a_i)$.
Otherwise let $G_i$ be an open neighbourhood of $W_{k-1}\backslash U^k_i$
in $W_{k-1} \cup U^k_i$ such that $\overline{G_i} \cap (U^k_i \backslash
W_{k-1})=\varnothing$. Let $F_i =\overline{G_i}\cap U^k_i$.

Observe that $U^k_i \cap W_{k-1}$ is $ANE([L])$ as an open subspace of
$AE([L])$-space $X$.
Hence $\cone (U^k_i \cap W_{k-1})$ is $AE([L])$ and therefore
inclusion of $F_i$ into the base of the cone can be extended to a map of
$U^k_i$ into this cone.
By $(ii)$ there exists an extension of $f_{k-1} |_{F_i}$ to the set $U^k_i$.
Let $f_k |_{U^k_i}$ be an extension of $f_{k-1} |_{F_i}$ such that
$\diam (f_k (U^k_i))<2\cdot\inf\{ \diam (g(U^k_i)) \mid g\mbox{ extends }
f_{k-1} |_{F_i}\}$.

Since $u_k$ is discrete system it suffices to check continuity of
$f_k$ at every point $a\in A$. Fix $\ep >0$.
Since $Y\in LC^{[L]}$ and $f_{k-1}$ is continuous mapping there exists
neighbourhood $O_a$ of $a$ in $X$ such that $f_{k-1} |_{O_a}$ is
$[L]$-homotopic to a constant map in $\ep/5$-neighbourhood of $f(a)$.
Applying property $(i)$ of $u$ to $O_a$ find neighbourhood $V_a$ of $a$.
Additionally, we may assume that $V_a =O(a,\delta)$ for some
$\delta >0$ such that $O(a,3\delta)\subset O_a$.
For all $U^k_i \in u_k$ such that $U^k_i \subset {\cal V}_k$
and $U^k_i \cap V_a \ne \varnothing$
we have $U^k_i \subset O_a$ by the choice of $V_a$.
Therefore construction of $f_k |_{U^k_i}$ and choice of $O_a$ imply
$\diam (f_k (U^k_i)) <\frac{4}{5} \ep$.
If $W_{k-1}\backslash U^k_i$ is open in $X$
then by the construction we have
$f(U^k_i)=f(a_i)$ where $a_i\in O_a$.
Hence $\dist (f(U^k_i),f(a))<\ep/5$ in this case.
Otherwise $f_k |_{U^k_i}$ was obtained as an extension of $f_{k-1}$
from nonempty set $F_i$ and it follows that
$\dist (f_k (U^k_i),f(a)) < \frac{4}{5}\ep + \frac{1}{5}\ep = \ep$.
Therefore $\dist (f_k (V_a),f(a))<\ep$ as required.
\end{proof}

The following theorem shows an importance of the notion of
lower $[L]$-continuity.
As an application of our selection theorem, we shall prove
the converse statement in section~\ref{S:appl}.

\begin{thm} \label{thmsoftness} Let $L$ be a $CW$-complex.
If a singlevalued continuous mapping $f\colon Y\to X$ of metric spaces
is locally $[L]$-soft, then the multivalued mapping
$f^{-1}\colon X\to Y$ is lower $[L]$-continuous.
If the mapping $f$ is $[L]$-soft, then
every fiber $f^{-1}(x)$ is $AE([L])$.
\end{thm}

\begin{proof}
Suppose that the mapping $f^{-1}\colon X\to Y$ is not lower $[L]$-continuous
at the point $\{x\}\times \{y\}$ of its graph.
Then there exist a positive $\ep$ and a sequence of mappings
$\{g_i\colon Z_i\to X,\ \wt g_i\colon A_i\to Y\}_{i=1}^\infty$,
where $A_i$ is a closed subset of paracompact
space $Z_i$ of extension dimension $\ed Z_i\le [L]$,
such that $f\circ \wt g_i=g_i|_{A_i}$, the images $g_i(Z_i)$
converges to the point $x$, the images $\wt g_i(A_i)$
converges to the point $y$, and the mapping $\wt g_i$ can not be
extended to a mapping of $Z_i$ into $O(y,\ep)$.

We consider a topological space $Z$ formed by the discrete union
of all spaces $Z_i$ and a point $\{p\}$ with the following topology:
an open base at the point $p$ consists of unions of this point and
all but finite number of spaces $Z_i$.
Clearly, the space $Z$ is paracompact and $\ed Z\le [L]$,
while the set $\{p\}\cup\bigcup_{i=1}^\infty A_i$ is closed in $Z$.
Let $g\colon Z\to X$ be a mapping such that $g|_{Z_i}=g_i$ and $g(p)=x$.
Also, let $\wt g\colon A\to Y$ be a mapping such that $g|_{A_i}=\wt g_i$
and $\wt g(p)=y$. These mappings are continuous and $f\circ \wt g=g|_A$.
It is easy to see that we can not extend the mapping $\wt g$
over neighbourhood of $A$ in $Z$ to a lifting of $g$ with respect to $f$.
Therefore, $f$ is not locally $[L]$-soft.
The first part of our lemma is proved.

Let the mapping $f\colon Y\to X$ be $[L]$-soft.
We consider a point $x$ and a mapping $h\colon A\to f^{-1}(x)$
of a closed subset $A$ of some paracompact space $Z$ with $\ed Z\le [L]$.
Since $f$ is $[L]$-soft, the constant mapping $h'\colon Z\to \{x\}$
admits a lifting $\wt h\colon Z\to f^{-1}(x)$ extending $h$.
Thus $f^{-1}(x)\in AE([L])$.
\end{proof}

\begin{thm}\label{thmmovement}
Let $L$ be a $CW$-complex such that $[L] \le [S^n]$ for some $n$.
Suppose that $F \colon X \to Y$ is a lower $[L]$-continuous multivalued
mapping of paracompact space $X$ to metric space $Y$. Let $K$ be a compact
subspace of a fiber $F(x)$ for some point $x \in X$. Then for any $\ep
>0$ there exist $\delta > 0$ and open neighbourhood $Ox$ of the point
$x$ such that for each $x' \in Ox$, for any paracompact space $Z$
with $\ed X \le [L]$, for each closed subspace $A$ of $Z$ and for any
map $f \colon (A,Z) \to (O(K,\delta) \bigcap F(x'), O(K,\delta))$
there exists $g \colon Z \to F(x')\bigcap O(K,\ep)$
such that $f|_A = g|_A$ and $\dist (f,g) < \ep$.
\end{thm}

\begin{proof}
Consider $\ep >0$. Using Lemma~\ref{lemmaequicompactset} choose sequence
$\{\delta_{-1}<\delta_0 <\delta_1 <\dots <\delta_n <\delta_{n+1} = \ep \}$
of positive numbers and neighbourhoods $\{O_i x \}_{i=0}^n$ of $x$
such that for all $i =-1, 0, 1, \dots, n$ and
for any points $x' \in O_i x$ and $y' \in F(x') \bigcap O(K,
\delta_i)$ the pair $O (y', \delta_{i}) \bigcap F(x') \subset
O(y', \delta_{i+1}/10) \bigcap F(x')$ is $[L]$-connected.
Let $\{ O(p_i, \delta_0/10) \mid i= 1, \dots,m \}$ be a finite
covering of compactum $K$ such that $p_i\in K$ for all $i$
and choose $\delta$ such that
$O(K,\delta)\subset \bigcup\limits_{i=1}^m O (p_i,\delta_0/10)$.
Let $Ox = \bigcap\limits_{i=1}^n O_i x$.

Fix $x' \in Ox$ and consider $f \colon Z \to O(K,\delta)$
such that $f(A) \subset F(x')\bigcap O(K,\delta)$ where $Z$ has
extension dimension $\ed Z \le [L]$.
Let $v$ be an open covering
$\{ V_p = f^{-1} O(p,\delta_0/10) \mid p = p_1, \dots, p_m \}$ of $Z$.
Find an open locally finite covering $\Sigma$ of $Z$
such that closures of elements of $\Sigma$ form strong star-refinement
of $v$ and order of $\Sigma$ is $\le n+1$.
For each $s \in \Sigma$ find $p(s) \in \{ p_1, \dots,p_m \}$
such that $\St (s,\Sigma) \subset V_{p(s)}\in v$
and pick $y_s \in O(p(s),\delta_0/10) \bigcap F(x')$.
Note that $f(s)\subset O(p(s),\delta_0/10)$.
Letting $g_{-1} = f|_A$ we shall inductively construct a sequence of mappings
$\{ g_k\colon\Sigma^{(k)}\cup A\to F(x')\}_{k=-1}^n$,
where $\Sigma ^{(k)}$ was defined in the beginning of Section 2,
such that $g_k$ extends $g_{k-1}$ and
$$g_k ((\Sigma^{(k)}\cup A)\cap s) \subset O(y_s,\delta_{k+1}/2)
\mbox{ for each } s\in\Sigma\eqno{(*)}$$
Since $\Sigma ^{(n)} =Z$ and $\delta_{n+1} =\ep$,
$(*)$ implies
$g_n (Z)\subset O(K,\delta_0/10+\ep
/2)\subset O(K,\ep)$.
Moreover,
$g_n$ is $\ep$-close to $f$, since for any $s\in\Sigma$ we have
$\dist (f|_s,g_n |_s) <
\dist (f|_s, p(s)) + \dist (p(s),y_s) +\dist (g_n |_s, y(s))
< \delta_0/10+ \delta_0/10
+ \ep/2 <\ep$.
Therefore, letting $g=g_n$ we shall obtain desired mapping.

Suppose that $g_k$ has been already constructed. It suffices to define
$g_{k+1}$ on the "interior" $\langle\sigma\rangle$ of each
"simplex" $[\sigma ] =[s_0, s_1, \dots, s_{k+1} ]$.
Let $[\sigma ]'=[\sigma ]\cap (\Sigma^{(k)}\cup A)$.
By property $(*)$ of $g_k$ we have
$\dist (g_k([\sigma]'), y_{s_0}) < \delta_{k+1}/2
+\max\limits_{i=1}^{k+1}\{\dist (y_{s_0}, y_{s_i})\}$.
Further, since $f(s)\subset O(p(s),\delta_0/10)$ for any $S$
and $s_0\cap s_i\ne \varnothing$, we have $\dist (p(s_0), p(s_i))
<2\delta_0/10$. Since
$y_{s_i}\in O(p(s_i),\delta_0/10)$, we therefore obtain

$\max\limits_{i=1}^{k+1}\{\dist (y_{s_0}, y_{s_i}) \} \le
\dist (y_{s_0}, p(s_0)) + \dist (p(s_0), p(s_i))
+\dist (p(s_i), y_{s_i}) <
\delta_0/10 + 2\delta_0/10 +\delta_0/10 = 2\delta_0/5$.
Therefore
$$g_k ([\sigma ]') \subset O(y_{s_0},\delta_{k+1}/2
+2\delta_0/5 ) \bigcap F(x')\subset O(y_{s_0},\delta_{k+1}
\bigcap F(x')$$
By the choice of $Ox$ and $\delta_{k+2}$
the pair

\[ O(y_{s_0},\delta_{k+1}) \bigcap F(x') \subset
O(y_{s_0},\delta_{k+2}/10)\bigcap F(x') \]

\noindent is $[L]$-connected.
Hence the map $g_k$ can be extended to a map
$g_{k+1}$ such that
$g_{k+1} ([\sigma ])\subset O(y_{s_0},\delta_{k+2}/10)\bigcap F(x')$.
Let us check the property $(*)$.
For any point $x\in (\Sigma^{(k)} \cup A)
\cap s_i$ by the construction of $g_{k+1}$ we have:
$\dist(g_{k+1}(x),y_{s_i}) <
\dist (g_{k+1} (x),y_{s_0}) +  \dist(y_{s_0}, y_{s_i})
\le \delta_{k+2}/10
+ 2(\delta_0/5) < \delta_{k+2}/2$, as required.
\end{proof}

\begin{cor} \label{cormovement}
Let $L$ be a $CW$-complex such that $[L] \le [S^n]$ for some $n$.
Let $Y$ be a metric space, $B$ be an $ANE([L])$-subspace of $Y$ and $K$
be a compact subspace of $B$. Then for any open neighbourhood $U$ of $K$
in $Y$ and for any $\ep > 0$ there exists a neighbourhood
$V \subset O(K,\ep)$ of $K$ with the following property:
for any paracompact space $X$ with $\ed X \le [L]$, any closed
subspace $A$ of $X$ and for any map $f \colon X \to V$ with $f(A) \subset B$
there exists a map $g \colon X \to U \cap B$ such
that $g$ is $\ep$-close to $f$ and $g|_A = f|_A$.
\end{cor}

\begin{lem} \label{lemmaequiUVcompactum}
Let $L$ be a $CW$-complex such that $[L]\le [S^n]$ for some $n$.
Let $F\colon X\to Y$ be lower $[L]$-continuous multivalued mapping
of topological space $X$ to metric space $Y$.
Suppose that a fiber $F(x)$ contains compact $UV^{[L]}$-pair $K \subset M$.
Then for any neighbourhood $U$ of $M$ in $Y$ there exist
neighbourhoods $V$ of $K$ in $Y$ and $O_x$ of the point $x$
in $X$ such that for any point $x'\in O_x$ the pair
$V \cap F(x')\subset U \cap F(x')$ is $[L]$-connected.
\end{lem}

\begin{proof}
Embed $Y$ into Banach space $E$ and consider $F$ as a mapping into $E$.
Fix $\ep>0$ and take a neighbourhood $O(M,3\ep)$ of $M$ in $E$.
By Theorem~\ref{thmmovement} there exist $\delta<\ep$ and a
neighbourhood $O_x$ of the point $x$ such that for any
point $x'\in O_x$, for any space $Z$ of extension dimension
$\ed Z\le [L]$ and its closed subset $A\subset Z$, and for any mapping
$\psi\colon (A,Z)\to (O(M,\delta)\cap F(x'),O(M,\delta))$
there exists a mapping $\psi'\colon Z\to F(x')$
such that $\psi'|_A=\psi|_A$ and $\dist(\psi,\psi')<\ep$.

Applying Homotopy Extension Theorem (see for example~\cite{Bor})
to $E$, we find a number $\sigma$ such that for any space $Z$,
any closed subspace $A$ of $Z$, and any two $\sigma$-close
maps $f,g\colon A\to O(K,\sigma)$ such that $f$ has an extension
$f'\colon Z\to O(M,\delta)$, it follows that $g$ also has
an extension $g'\colon Z\to O(M,2\delta)$ which is $\delta$-close to $f'$.
Using the $UV^{[L]}$-property of the pair $K\subset M$ in $F(x)$,
we take a number $\mu<\sigma$ such that the pair
$O(K,\mu)\cap F(x)\subset O(K,\delta)\cap F(x)$ is $[L]$-connected.
By Theorem~\ref{thmmovement} there exists $\nu<\mu$ such that for any
space $A$ of extension dimension $\ed A\le [L]$ and for any mapping
$\varphi\colon A\to O(K,\nu)$ there is a mapping
$\varphi'\colon A\to O(K,\mu)\cap F(x)$ with $\dist (\varphi,\varphi')<\mu$.
Put $V=O(K,\nu)$.

Consider a point $x'\in O_x$, a space $Z$ of extension dimension
$\ed Z\le [L]$ and its closed subspace $A\subset Z$.
Now any mapping $\varphi\colon A\to V\cap F(x')$ is $\mu$-close
to some mapping $\varphi'\colon A\to O(K,\mu)\cap F(x)$
which can be extended to a mapping
$\wt\varphi'\colon Z\to O(M,\delta)\cap F(x)$.
Since $\varphi|_A$ and $\varphi'|_A$ are $\sigma$-close maps into
$O(K,\sigma)$, $\varphi$ can also be extended to a mapping
$\psi\colon Z\to O(M,2\delta)$ which is $\delta$-close to $\wt\varphi'$.
Finally, there is another extension
$\psi'\colon Z\to O(M,2\delta+\ep)\cap F(x')$ of the mapping $\varphi$.
Thus, the pair $V\cap F(x')\subset O(M,3\ep)\cap F(x')$ is $[L]$-connected.
\end{proof}

\begin{lem} \label{lemUVLpair}
Let $L$ be a $CW$-complex such that $[L]\le [S^n]$ for some $n$.
Consider spaces $K \subset M \subset Y \subset E$, where $K$ and $M$
are compacta, $Y$ and $E$ are metric $ANE ([L])$-spaces. Then $K\subset M$ is
$UV^{[L]}$-pair in $Y$ if and only if it is $UV^{[L]}$-pair in $E$.
\end{lem}

\begin{proof}
If $K\subset M$ is $UV^{[L]}$-pair in $Y$, consider a multivalued mapping
$F$ of the unit interval $I=[0,1]$ defined as follows:
$F(0)=Y$ and $F(x)=E$ for any positive $x\in I$.
Clearly, $F$ is lower $[L]$-continuous.
Now Lemma~\ref{lemmaequiUVcompactum} implies the $UV^{[L]}$-property of
the pair $K\subset M$ in $E$.

Assume that $K\subset M$ is $UV^{[L]}$-pair in $E$.
Take an open neighbourhood $U$ of $M$ in $Y$ and consider
an open neighbourhood $O(M,2\ep)$ in $E$ such that
$O(M,2\ep)\cap Y\subset U$.
By Corollary~\ref{cormovement} there exists $\delta<\ep$
such that for any space $Z$ of extension dimension
$\ed Z\le [L]$ and its closed subset $A\subset Z$, and for any mapping
$\psi\colon (A,Z)\to (O(K,\delta)\cap Y,O(K,\delta))$
there exists a mapping $\psi'\colon Z\to Y$
such that $\psi'|_A=\psi|_A$ and $\dist(\psi,\psi')<\ep$.
Using the $UV^{[L]}$-property of the pair $K\subset M$ in $E$,
we can find a neighbourhood $V'$ of $K$ in $E$.
Put $V=V'\cap Y$.

Now any mapping $\varphi \colon A\to V$ of closed subset $A$
of space $Z$ of extension dimension $\ed Z\le [L]$ can
be extended to a mapping $\psi \colon Z\to O(K,\delta)$.
And by the choice of $\delta$ there is an extension
$\psi'\colon Z\to O(M,2\ep)\cap Y$ of the mapping $\varphi$.
\end{proof}

\begin{thm} \label{lemmainvar}
Let $L$ be a $CW$-complex such that $[L]\le [S^n]$ for some $n$.
Suppose that a compact pair $K \subset M$ is $UV^{[L]}$-connected
with respect to embedding in some Polish $ANE([L])$-space $B$.
Then this pair is $UV^{[L]}$-connected with respect to any embedding
in any Polish $ANE([L])$-space.
\end{thm}

\begin{proof}
There exists an embedding $i \colon M \to \R^{\omega}$ which can be
extended to an embedding of any Polish space containing $M$
(see Theorem~2.3.17 in~\cite{Ch96}).

If the pair $K \subset M$ is $UV^{[L]}$-connected in a Polish space $B$,
then we can extend $i$ to an embedding of $B$ in $\R^{\omega}$ and the
pair $K\subset M$ is $UV^{[L]}$-connected in $\R^{\omega}$
by Lemma~\ref{lemUVLpair}.

Consider any Polish $ANE([L])$-space $Y$, containing $M$. Extending $i$
to an embedding of $Y$ into $\R^{\omega}$, we obtain
$UV^{[L]}$-connectedness of the pair $K\subset M$ in $Y$ by
Lemma~\ref{lemUVLpair}.
\end{proof}

\section{Compact-valued selections}\label{S:compactselections}

This section is devoted to the construction of compact-valued upper
semicontinuous selections for multivalued mappings.

\begin{lem} \label{lemmapairs}
Let $f\colon X\to Y$ be a continuous singlevalued
mapping of compact metric spaces.
Let $Y_1\subset Y$ be a closed subset and $X_1$ be its inverse
image $X_1=f^{-1}(Y_1)$.
If the mapping $f|_{X_1}\colon X_1\to Y_1$ is approximately
$[L]$-invertible and the pair $X_1\subset X$ is $UV^{[L]}$-connected,
then the pair $Y_1\subset Y$ is also $UV^{[L]}$-connected.
\end{lem}

\begin{proof}
Consider $f$ as a submapping of
the projection $\pi\colon l_2\times l_2\to l_2$.
Let $U$ be some neighbourhood of a compact space $Y$ in $l_2$.
We must find a neighbourhood $V$ for $Y_1$ such that the pair
$V\subset U$ is $[L]$-connected.

By the $UV^{[L]}$-connectedness of the pair $X_1\subset X$, we fix
an open neighbourhood $W$ of $X_1$ such that the pair
$W\subset \pi^{-1}(U)$ is $[L]$-connected.
By approximate $[L]$-invertibility of the mapping $f|_{X_1}$
there exists a neighbourhood $V$ of $Y_1$ such that any mapping
$g\colon Z\to V$ of the space $Z$ of extension dimension
$\ed Z\le [L]$ admits a lifting map $\wt g\colon Z\to U$.

Now if $g\colon A\to V$ is a mapping of closed subset $A\subset Z$
where $\ed Z\le [L]$, we take a lifting map
$\wt g\colon A\to W$ and extend it to a mapping $g'\colon Z\to \pi^{-1}(U)$.
Define an extension of $g$ as $\pi\circ g'$.
\end{proof}

By $\exp Z$ is denoted the space of all compact subsets of a metric
space $Z$ endowed with the Hausdorff metric.

\begin{defin}
The {\it exponential of a pair} $\exp (A,B)$ is a subspace of $\exp B$
formed by compact sets $K\subset B$ containing $A$.
We define the {\it $UV^{[L]}$-exponential of the pair} $(A,B)$
as follows:
$$ UV^{[L]}\mbox{-}\exp(A,B)=\{ K\in \exp B \mid
  \mbox{ the pair } A\subset K \mbox{ is $UV^{[L]}$-connected}\}. $$
\end{defin}

\begin{lem} \label{lemmaclosedness}
For any pair $(K,X)$ formed by a compact set $K$ and a metric space $X$,
the set ${\rm UV}^{[L]}\mbox{-}\exp (K,X)$ is closed in $\exp (K,X)$.
\end{lem}

\begin{proof}
Let a sequence of compact sets $\{K_m\}_{m\ge 1}$ from the
${\rm UV}^{[L]}$-exponential of the pair $(K,X)$ be convergent with
respect to the Hausdorff metric to a compact set $K_0$.
Consider a neighbourhood $U$ of $K_0$.
There exists $m\ge 1$ such that $K_m\subset U$.
Now $UV^{[L]}$-connectedness of the pair $K\subset K_m$
allows us to find a neighbourhood $V$ of the compact
set $K$ such that the pair $V\subset U$ is $[L]$-connected.
\end{proof}

\begin{defin}
The {\it fiberwise exponential} of a multivalued
mapping $F\colon X\to Y$ is the mapping $\exp F\colon X\to \exp Y$
which assigns $\exp F(x)$ to a point $x$.
\end{defin}

\begin{lem} \label{lemmacompleteness}
The fiberwise exponential of a complete mapping is complete.
\end{lem}

\begin{proof}
Since the exponential of an open set is open and the
exponential of an intersection coincides with the intersection of
exponentials, the exponential of a $G_\delta$-set is a
$G_\delta$-set. Since the exponential of a closed set is closed,
the exponential of a fiber closed in a $G_\delta$-set is closed in
the exponential of a $G_\delta$-set.
\end{proof}

\begin{lem} \label{lemmacompactsupport}
Let $L$ be a finite $CW$-complex
such that $[L]\le [S^n]$ for some $n$.
Suppose that a metric space $Z$ contains a compactum $K$
and the pair $K\subset Z$ is $[L]$-connected.
Then there exists a compactum $K'\subset Z$ containing $K$
such that the pair $K\subset K'$ is $UV^{[L]}$-connected.
\end{lem}

\begin{proof}
By proposition 2.23 in \cite{Ch},
there is a compactum $X$ of extension dimension $\ed X\le [L]$
and a continuous mapping $f$ of $X$ onto $K$ such that every fiber
$f^{-1} (y)$ is $UV^{[L]}$--compactum.
By Theorem \ref{thmapprinvert},
the mapping $f$ is approximately $[L]$--invertible.
There exists $AE([L])$-compactum $X^{[L]}$ containing $X$ such
that $\ed X^{[L]}=[L]$~\cite{Ch97}.
It is easy to see from Lemma~\ref{lemmainvar} that
the pair $X\subset X^{[L]}$ is $UV^{[L]}$-connected.

Since the pair $K\subset Z$ is $[L]$-connected, we can extend the mapping
$f$ to a mapping $\wt f\colon X^{[L]}\to Z$. Put $K'=\wt f(X^{[L]})$.
Then the pair $K\subset K'$ is $UV^{[L]}$-connected by Lemma~\ref{lemmapairs}.
\end{proof}

\begin{defin}
For a multivalued mapping $\Phi\colon X\to Y$ and its compact
submapping $\Psi$ we define {\it fiberwise $UV^{[L]}$-exponential of the pair}
$UV^{[L]}\mbox{-}\exp(\Psi,\Phi)\colon X\to \exp Y$ as a mapping
assigning $UV^{[L]}\mbox{-}\exp(\Psi(x),\Phi(x))$ to a point $x\in X$.
\end{defin}

\begin{lem} \label{lemmacontinuity}
Let $L$ be a finite $CW$-complex
such that $[L]\le [S^n]$ for some $n$.
Suppose that a lower $[L]$-continuous mapping $\Phi\colon X\to Y$
of paracompact space $X$ to metric space $Y$
contains a compact submapping $\Psi$.
Then the fiberwise ${\rm UV}^{[L]}$-exponential of the pair
${\rm UV}^{[L]}\mbox{-}\exp(\Psi,\Phi)$ is lower semicontinuous.
\end{lem}

\begin{proof}
Denote $F={\rm UV}^{[L]}\mbox{-}\exp (\Psi,\Phi)$, and
for a point $x\in X$ fix a compact set $K\in F(x)$.
Fix a positive number $\ep$.
By Lemma~\ref{lemmaequiUVcompactum} there are number $\delta<\ep$
and neighbourhood $O'_x$ of the point $x$ such that the pair
$O(\Psi(x),\delta) \cap \Phi(x')\subset O(K,\ep)\cap \Phi(x')$
is $[L]$-connected for any point $x'\in O'_x$.
Since $\Phi$ is lower semicontinuous and $K$ is compact, there
exists a neighbourhood $O''_x$ of the point $x$ such that
$O(y,\ep/2)\cap \Phi(x')\ne\emptyset$ for any points
$y\in K$ and $x'\in O''_x$ (apply Lemma~\ref{lemmaequicompactset}).
Let $O_x$ be a neighbourhood of $x$ such that
$O_x\subset O'_x\cap O''_x$ and $\Psi(x')\subset O(\Psi(x),\delta)$
for every point $x'\in O_x$.

Take any point $x'\in O_x$. By Lemma~\ref{lemmacompactsupport}
there exists a compactum $\wt K\subset \Phi(x')\cap O(K,\ep)$
such that the pair $\Psi(x')\subset \wt K$ is $UV^{[L]}$-connected,
and therefore $\wt K\in F(x')$.
It remains to enlarge (if necessary) the compactum $\wt K$
to obtain a compactum $K'$ with $\dist(\wt K,K')<\ep$.
By the choice of the neighbourhood $O''_x$ there is a finite
set of points $P$ in $\Phi(x')$ such that $\dist (K,P)<\ep$.
We put $K'=\wt K\cup P$.
\end{proof}

\begin{lem} \label{lemmafiberwisecontraction}
Let $L$ be a finite $CW$-complex
such that $[L]\le [S^n]$ for some $n$.
Let $\Phi\colon X\to Y$ be a complete lower $[L]$-continuous mapping
of a paracompact space into a complete metric space containing a compact
submapping $\Psi$ such that the pair $\Psi\subset\Phi$
is fiberwise $[L]$-connected.
Then there exists a compact submapping $\Psi'$ of the mapping
$\Phi$ such that the pair $\Psi(x)\subset \Psi'(x)$
is $UV^{[L]}$-connected for any $x\in X$.
\end{lem}

\begin{proof}
Consider $F={\rm UV}^{[L]}\mbox{-}\exp (\Psi,\Phi)$.
According to Lemma~\ref{lemmacompactsupport}, the
mapping $F$ has nonempty fibers.
By Lemma~\ref{lemmacontinuity}, $F$ is lower semicontinuous.
By Lemma~\ref{lemmaclosedness}, $F$ is fiberwise closed in
$\exp (\Psi,\Phi)$, and therefore, the completeness
of this mapping follows from the completeness of the latter,
which was established in Lemma~\ref{lemmacompleteness}.
Then by the compact-valued selection theorem from~\cite{SB},
the mapping $F$ admits a compact selection $F'$.
Define a compact mapping $\Psi'\colon X\to Y$ by the
equality $\Psi'(x)=\bigcup_{K\in F'(x)} K$. Since for any $K\in F'(x)$,
the pair $\Psi(x)\subset K$ is $UV^{[L]}$-connected, then
the pair $\Psi(x)\subset\Psi'(x)$ is also $UV^{[L]}$-connected.
\end{proof}

\begin{lem} \label{lemmacompactfiltration}
Let $L$ be a finite $CW$-complex
such that $[L]\le [S^n]$ for some $n$.
Then any $[L]$-connected lower $[L]$-continuous increasing
$n$-filtration $\Phi=\{\Phi_k\}$ of complete mappings of a
paracompact space to a complete metric space contains a compact
${\rm UV}^{[L]}$-connected $n$-subfiltration $\Psi=\{\Psi_k\}$.
\end{lem}

\begin{proof}
The construction of filtration $\Psi$ is performed by
induction with the use of Lemma~\ref{lemmafiberwisecontraction}.
The initial step of induction consists in the construction
of a compact submapping $\Psi_0\subset\Phi_0$.
This can be done by the use of the compact-valued selection theorem
from~\cite{SB} since
the initial term of the filtration $\Phi$ is lower semicontinuous.
If compact ${\rm UV}^{[L]}$-connected filtration
$\{\Psi_m\}_{m<k}$ has been constructed such that $\Psi_m\subset \Phi_m$
for $m<k$, then the pair $\Psi_{k-1}\subset\Phi_{k}$ satisfies the
conditions of Lemma~\ref{lemmafiberwisecontraction}, and according to this
lemma, we complete the construction of the filtration.
\end{proof}

The following lemma is a generalization of Lemma~\ref{lemmaequiUVcompactum}
and we will use it in section~\ref{S:appl}.

\begin{lem} \label{equiUVmap}
Let $L$ be a $CW$-complex such that $[L]\le [S^n]$ for some $n$.
Let $F\colon X\to Y$ be lower $[L]$-continuous multivalued mapping
of paracompact space $X$ to metric space $Y$.
For a closed subset $A\subset X$ consider a compact submappings
$H\subset \wt H\colon A\to Y$ of the mapping $F|_A$.
If the pair $H\subset \wt H$ is fiberwise $UV^{[L]}$-connected, then for any
neighbourhood $\mathcal U$ of the graph $\Gamma_{\wt H}$ in the product
$X\times Y$ there exists a neighbourhood $\mathcal V$ of the graph
$\Gamma_{H}$ in the product $X\times Y$ such that the pair
$\mathcal V(x)\cap F(x)\subset \mathcal U(x)\cap F(x)$ is
$[L]$-connected for every $x$ from some open neighbourhood of the set $A$.
\end{lem}

\begin{proof}
By Lemma~\ref{lemmaequiUVcompactum} we take for every point $x\in A$
an open neighbourhood $O_x\subset X$ of the point $x$ and
an open neighbourhood $V_x\subset Y$ of the set $H(x)$
such that the set $H(O_x\cap A)$ is contained in $V_x$
and the pair $V_x\cap F(x')\subset \mathcal U(x')\cap F(x')$
is $[L]$-connected for every point $x'\in O_x$.
Fix a closed neighbourhood $B$ of the set $A$ such that
$B\subset \cup_{x\in A} O_x$.
Let $\Omega_1=\{\omega_\lambda\}_{\lambda\in\Lambda}$
be a locally finite open (in $B$) cover of $B$ refining
the cover $\{O_x\}_{x\in A}$.
For every $\lambda\in \Lambda$ we take a set $V_\lambda=V_x$
such that $\omega_\lambda\subset O_x$.
Let $\Omega_2\in \cov B$ be a locally finite open cover starlike
refining $\Omega_1$. For $x\in\Int B$ we define
$$\mathcal V(x)=\cap\{V_\lambda\mid \St(x,\Omega_2)\subset\omega_\lambda\}. $$
Since the cover $\Omega_1$ is locally finite, the set $\mathcal V(x)$ is an
intersection of finitely many open sets, and, therefore,
$\mathcal V(x)$ is open.

Since for every $\lambda$ the pair
$V_\lambda\cap F(x)\subset \mathcal U(x)\cap F(x)$ is $[L]$-connected, then the
pair $\mathcal V(x)\cap F(x)\subset \mathcal U(x)\cap F(x)$ is $[L]$-connected.
Since the cover $\Omega_2$ is locally finite, then for every point
$x\in\Int B$ there is a neighbourhood $W_x$ such that for
any point $x'\in W_x$ we have $\St(x,\Omega_2)\subset\St(x',\Omega_2)$.
Therefore, for every $x'\in W_x$ we have $\mathcal V(x)\subset \mathcal V(x')$.
Thus, the set $\mathcal V$ is open.
\end{proof}

\begin{cor} \label{cornhdselection}
Let $L$ be a $CW$-complex such that $[L]\le [S^n]$ for some $n$.
Suppose that lower $[L]$-continuous multivalued mapping $F\colon X\to Y$
of paracompact space $X$ to metric space $Y$
contains a singlevalued continuous selection $f\colon A\to Y$
over the closed subset $A\subset X$.
Then for any neighbourhood $\mathcal U$ of the graph $\Gamma_f$ in the product
$X\times Y$ there exists a neighbourhood $\mathcal V$ of the graph $\Gamma_f$
in the product $X\times Y$ such that for every point $x\in \pr_X \mathcal V$
the pair $\mathcal V(x)\cap F(x)\subset \mathcal U(x)\cap F(x)$
is $[L]$-connected.
\end{cor}

\section{Selection theorems}\label{S:paracom}

The {\it gauge} of a multivalued mapping $F\colon X\to Y$ is defined as
 $$\operatorname{cal} (F)=\sup\{\operatorname{diam} F(x) |\; x\in X\}.$$

\begin{lem} \label{lemmashrinking}
Let $L$ be a finite $CW$-complex
such that $[L]\le [S^n]$ for some $n$.
Let $X$ be a paracompact space of extension dimension $\ed X\le [L]$.
If a complete lower $[L]$-continuous mapping $\Phi\colon X\to Y$
into a complete metric space $Y$ contains an $n$-$UV^{[L]}$-filtered compact
submapping $\Psi$, then any neighbourhood of
the graph $\Gamma_\Psi$ contains the graph of a compact $n$-$UV^{[L]}$-filtered
submapping $\Psi'$ of the mapping $\Phi$ whose gauge
$\calibr(\Psi')$ does not exceed any given $\varepsilon$.
\end{lem}

\begin{proof}
Given an arbitrary number $\varepsilon >0$ and an open neighbourhood
$\mathcal U$ of the graph $\Gamma_\Psi$ in the product $X\times Y$,
consider a covering $\omega_n\times\varepsilon_n$ of
the graph $\Gamma_\Psi$ such that the star
$\St(\Gamma_\Psi,\omega_n\times\varepsilon_n)$ is contained in $\mathcal U$
(Lemma~\ref{lemmarectcover} is applied), while the function
$\varepsilon_n(x)$ does not exceed $\varepsilon/3$.

For an $[L]$-continuous mapping $\Phi$ and for its compact submapping $\Psi$,
applying successively Lemma~\ref{lemmastarlikerefinement}, we construct the
coverings $\{\omega_k\times\varepsilon_k\}_{k=0}^{n-1}$ such that
$\omega_k\times\varepsilon_k$ is starlike $[L]$-connectedly refined
into $\omega_{k+1}\times\varepsilon_{k+1}$ for any $k<n$.
By Lemma~\ref{lemmaequicompactmap} there is a neighbourhood $\mathcal V$ of
the graph $\Gamma_\Psi$ in the product $X\times Y$ such that
for any point $\{x\}\times\{y\}\in \mathcal V$, the star of this point
relative to the covering $\omega_0\times\varepsilon_0$ intersects the
fiber $\{x\}\times\Phi(x)$.

By Theorem~\ref{thmapproximation}, there is a continuous singlevalued
mapping $\psi\colon X\to Y$ whose graph is contained in $\mathcal V$.
We fix an $[L]$-connected $n$-filtration $\{G_m\}$ given
fiberwise by the equality
$G_m(x)=\Phi(x)\cap \St(\{x\}\times\psi(x),\omega_m\times\varepsilon_m)(x)$.
Since the projection of the star
$\St(\{x\}\times\psi(x),\omega_n\times\varepsilon_n)$ onto $Y$ has the
diameter less than $\varepsilon$, then $\calibr G_n<\varepsilon$.
By Lemma~\ref{lemmastartrace} the filtration
$\{G_m\}$ is complete and lower $[L]$-continuous.
Finally, Lemma~\ref{lemmacompactfiltration} allows us to find a compact
$UV^{[L]}$-connected $n$-subfiltration $\Psi'=\{\Psi'_k\}$.
\end{proof}

\begin{thm} \label{thmmixedselection}
Let $L$ be a finite $CW$-complex
such that $[L]\le [S^n]$ for some $n$.
Let $X$ be a paracompact space of extension dimension $\ed X\le [L]$.
If a complete lower $[L]$-continuous multivalued mapping $\Phi$
of $X$ into a complete metric space $Y$ contains an
$n$-$UV^{[L]}$-filtered compact submapping $\Psi$,
then $\Phi$ contains a singlevalued continuous selection $s$.

A selection $s$ can be chosen in such a way that the graph of this
selection is contained in any given neighbourhood $\mathcal U$ of the graph
$\Gamma_\Psi$ in the product $X\times Y$.
\end{thm}

\begin{proof}
Let $\{\Psi_k\}_{k=0}^n$ be $UV^{[L]}$-filtration of $\Psi$.
Denote $\Psi_n$ by $\Psi_n^0$ and take an arbitrary neighbourhood
$\mathcal U_0$ of the graph $\Psi_n^0$.
Consider a $G_\delta$-subset $G\subset X\times Y$ such that
all fibers of $F$ are closed in $G$ and fix open sets
$G_i\subset X\times Y$ such that $G=\cap^\infty_{i=0} G_i$.
By induction with the use of Lemma~\ref{lemmashrinking},
we construct a sequence of $n$-$UV^{[L]}$-filtered mappings
$\{\Psi_n^k\}_{k=1}^\infty$ and of open neighbourhoods of graphs
of these mappings $\{\mathcal U_k\}_{k=1}^\infty$ such that for any $k\ge 1$,
the gauge $\calibr\Psi_n^k$ does not exceed $\frac{1}{2^k}$, while
the graph $\Psi^k_n$ together with its neighbourhood
$\mathcal U_k$ is in $\mathcal U_{k-1}\cap G_{k-1}$.
It is not difficult to choose the neighbourhood $\mathcal U_k$ of the graph
$\Psi_n^k$ in such a way that the fibers $\mathcal U_k(x)$
have the diameter not more than $3\cdot\calibr\Psi_n^k=\frac{3}{2^k}$.

Then for any $m\ge k\ge 1$ and for any point $x\in X$,
$\Psi_n^m(x)\subset O(\Psi_n^k(x),\frac{3}{2^k})$; this implies
that $\{\Psi_n^k\}_{k=1}^\infty$ is a Cauchy sequence.
Since $Y$ is complete, there exists a limit $s$ of this sequence.
The mapping $s$ is singlevalued by
the condition $\calibr\Psi_n^k<\frac{1}{2^k}$ and is upper
semicontinuous (and, therefore, is continuous) by the upper
semicontinuity of all the mappings $\Psi_n^k$.
Clearly, for any $x\in X$ the point $s(x)$ lies in $G(x)$
and is a limit point of the set $F(x)$.
Since $F(x)$ is closed in $G(x)$, then $s(x)\in F(x)$,
i.e. $s$ is a selection of the mapping $F$.
\end{proof}

\begin{cor}
Let $L$ be a finite $CW$-complex
such that $[L]\le [S^n]$ for some $n$.
Let $X$ be a paracompact space of extension dimension $\ed X\le [L]$.
Let a complete lower $[L]$-continuous multivalued mapping $\Phi$ of $X$ into
a complete metric space $Y$ contains an $n$-$UV^{[L]}$-filtered compact
submapping $\Psi$ which is singlevalued on some closed subset $A\subset X$.
Then any neighbourhood $\mathcal U$ of the graph $\Gamma_\Psi$ in the product
$X\times Y$ contains the graph of a singlevalued continuous selection
$s$ of the mapping $\Phi$ which coincides with $\Psi\big|_A$ on the set $A$.
\end{cor}

\begin{proof}
Apply Theorem~\ref{thmmixedselection} to the mapping $F$ defined as follows:
$$
  F(x)=\cases \Psi(x),
         &\text{if \;$x\in A$}\\
        \Phi(x),
         &\text{if \;$x\in X\setminus A.$}\endcases
$$
\end{proof}

\begin{thm} \label{thmfilteredselection}
Let $L$ be a finite $CW$-complex
such that $[L]\le [S^n]$ for some $n$.
Let $X$ be a paracompact space of extension dimension $\ed X\le [L]$.
Suppose that multivalued mapping $F\colon X\to Y$
into a complete metric space $Y$ admits a lower $[L]$-continuous, complete,
and fiberwise $[L]$-connected $n$-filtration
$F_0\subset F_1\subset\dots\subset F_n\subset F$.
If $f\colon A\to Y$ is a continuous singlevalued selection of $F_0$
over a closed subspace $A\subset X$, then there exists a continuous
singlevalued selection $\wt f\colon X\to Y$
of the mapping $F$ such that $\wt f|_A=f$.
\end{thm}

\begin{proof}
For every $i\le n$ define a multivalued mapping $\Phi_i\colon X\to Y$
as follows:
$$
  \Phi_i(x)=\cases f(x),
         &\text{if \;$x\in A$}\\
        F_i(x),
         &\text{if \;$x\in X\setminus A.$}\endcases
$$
Then $\Phi_0\subset\Phi_1\subset\dots\subset\Phi_n$ is lower
$[L]$-continuous, complete, and fiberwise $[L]$-connected $n$-filtration.
By Lemma~\ref{lemmacompactfiltration} the mapping $\Phi_n$ contains
a compact ${\rm UV}^{[L]}$-connected $n$-subfiltration.
And application of Theorem~\ref{thmmixedselection} completes the proof.
\end{proof}

\begin{thm} \label{thmlocalselection}
Let $L$ be a finite $CW$-complex
such that $[L]\le [S^n]$ for some $n$.
Let $X$ be a paracompact space of extension dimension $\ed X\le [L]$.
Let $F\colon X\to Y$ be a complete lower $[L]$-continuous multivalued
mapping into a complete metric space.
Suppose that $f\colon A\to Y$ is a continuous singlevalued selection of $F$
over a closed subspace $A\subset X$. Then there exists a continuous
extension of $f$ to a selection of the mapping $F$ over some neighbourhood of
the set $A$.
\end{thm}

\begin{proof}
Put $\mathcal U_n=X\times Y$. Using Corollary~\ref{cornhdselection}
we find open neighbourhoods
$\mathcal U_0\subset\mathcal U_1\subset\dots\subset\mathcal U_n$
of the graph $\Gamma_f$ in $X\times Y$ such that for any
$x\in \pr_X \mathcal U_0$ the pair
$\mathcal U_i(x)\cap F(x)\subset\mathcal U_{i+1}(x)\cap F(x)$
is $[L]$-connected for every $i<n$.
Let $OA$ be a closed neighbourhood of $A$ contained in $\pr_X \mathcal U_0$.
For every $i\le n$ define a multivalued mapping $F_i\colon OA\to Y$
by equality $F_i(x)=\mathcal U_i(x)\cap F(x)$.
Then $F_0\subset F_1\subset\dots\subset F_n=F|_{OA}$
is fiberwise $[L]$-connected $n$-filtration.
As a closed subset of $X$, the space $OA$ is paracompact
of extension dimension $\le [L]$.
It is easy to see that every mapping $F_i$ is lower
$[L]$-continuous and complete.
Applying Theorem~\ref{thmfilteredselection} we extend $f$
to a selection of $F$ over $OA$.
\end{proof}

\section{Applications of selection theorems}\label{S:appl}

The following theorem is well-known for $n$-soft mappings~\cite{Dr1}.

\begin{thm} \label{thmcharLsoftness}
Let $L$ be a finite $CW$-complex
such that $[L]\le [S^n]$ for some $n$.
A singlevalued continuous mapping $f\colon Y\to X$ of Polish
spaces is locally $[L]$-soft if and only if the multivalued mapping
$f^{-1}\colon X\to Y$ is lower $[L]$-continuous.
The mapping $f$ is $[L]$-soft if and only if
every fiber $f^{-1}(x)$ is $AE([L])$ and
the mapping $f^{-1}$ is lower $[L]$-continuous.
\end{thm}

\begin{proof}
The part "only if" is proved in section~\ref{S:soft}
(Theorem~\ref{thmsoftness}).

For the "if" part, consider a paracompact space $Z$ with
$\ed Z\le [L]$, its closed subset $A\subset Z$, and continuous mappings
$g\colon Z\to X$ and $g'\colon A\to Y$ such that $g|_A=f\circ g'$.
Then the multivalued mapping $F\colon Z\to Y$ defined as
$F=f^{-1}\circ g$ is lower $[L]$-continuous and complete.
By Theorem~\ref{thmlocalselection} a selection $g'$ of $F$
admits an extension $\wt g$ on some open neighbourhood $OA$ of the set $A$.
If every set $f^{-1}(x)$ is $AE([L])$, then filtration
$F\subset F\subset\dots\subset F$ is fiberwise $[L]$-connected
and by Theorem~\ref{thmfilteredselection} we can assume
that $\wt g$ is defined on $Z$. Clearly, $\wt g$ is a lifting
of $g$ and theorem is proved.
\end{proof}

\begin{lem} \label{lemmaextUVL}
Let $L$ be a finite $CW$-complex such that $[L]\le [S^n]$ for some $n$.
Let $F\colon X\to Y$ be lower $[L]$-continuous complete multivalued mapping
of a separable metric space $X$ with $\ed X\le [L]$ to Polish space $Y$.
Suppose that $\Psi\colon A\to Y$ is u.s.c. $UV^{[L]}$-valued submapping
of $F|_A$ defined on closed subset $A\subset X$.
Then there exists u.s.c. $UV^{[L]}$-valued submapping $\Psi'\colon OA\to Y$
of $F|_{OA}$ defined on some neighbourhood $OA$ of $A$ such that
$\Psi'|_A=\Psi$, and $\Psi'|_{OA\setminus A}$ is continuous and singlevalued.
\end{lem}

\begin{proof}
Using Lemma~\ref{equiUVmap}, we can construct a sequence
$\{\mathcal U_i\}^\infty_{i=1}$ of open in $X\times Y$
neighbourhoods of the graph $\Gamma_\Psi$
such that $\mathcal U_0=X\times Y$ and for every $i\ge 0$
the pair $\mathcal U_{i+1}(x)\cap F(x)\subset \mathcal U_{i}(x)\cap F(x)$
is $[L]$-connected for all
$x$ from some open neighbourhood $O_iA$ of the set $A$.
We may assume that the set $\mathcal U_i$ is contained in
${\frac{1}{2^{i+1}}}$-neighbourhood of the graph $\Gamma_\Psi$
(for metric spaces $(X,\rho_X)$ and $(Y,\rho_Y)$
we equip the product $X\times Y$ with a metric
$\rho((x_1,y_1),(x_2,y_2))=\max\{\rho_X(x_1,x_2),\rho_Y(y_1,y_2)\}$).

Take a sequence $\{F_k\}_{k=1}^\infty$ of closed neighbourhoods
of the set $A$ such that $F_k\subset pr_X(\mathcal U_k)\cap O_{k-1}A$ and
$F_{k+1}\subset\Int (F_k)$ for every $k\ge 1$. Put $OA=F_{n+1}$.
Define the maps $\{\Phi_m\colon F_{n}\setminus A\to Y\}^{n}_{m=0}$ by the rule
$\Phi_m(x)=\mathcal U_{k-m}(x)\cap F(x)$ for all $x\in F_k\setminus F_{k+1}$.
Using Theorem~\ref{thmfilteredselection},
we obtain a continuous singlevalued selection
$f\colon OA\setminus A\to Y$ of the map $\Phi_{n}$.
Let the map $\Psi'\colon OA\to Y$ be given by $\Psi$ on $A$ and
by $f$ on $OA\setminus A$.
Since the graph $\Gamma_f$ over the set
$F_k\setminus A$ is contained in $\mathcal U_{k-n-1}$
(and, therefore, in ${\frac{1}{2^{k-n}}}$-neighbourhood
of the graph $\Gamma_\Psi$),
we see that $\Psi'$ is upper semicontinuous.
\end{proof}

\begin{thm} \label{thmextUVL}
Let $L$ be a finite $CW$-complex such that $[L]\le [S^n]$ for some $n$.
Let $\Psi\colon A\to l_2$ be u.s.c. $UV^{[L]}$-valued mapping of a
closed subset $A\subset X$ of separable metric space $X$.
Then there exists u.s.c. $UV^{[L]}$-valued mapping $\Psi'\colon X\to l_2$
such that $\Psi'|_A=\Psi$.
\end{thm}

\begin{proof}
Consider proper continuous mapping $f\colon Y\to X$ of separable metric
spaces such that every fiber $f^{-1}(x)$ is $UV^{[L]}$-compactum
and $\ed Y\le [L]$ (see proposition 2.23 in~\cite{Ch}).
Denote by $A'$ the set $f^{-1}(A)\subset Y$.
Using Lemma~\ref{lemmaextUVL} we can find u.s.c. $UV^{[L]}$-valued
extension $F\colon Y\to l_2$ of the mapping $\Psi\circ f\colon A'\to l_2$
which is singlevalued and continuous on $Y\setminus A'$.
Let $\beta$ be positive continuous function on $Y\setminus A'$
such that $\beta(y)=\dist(f(y),A)$.
Using propositions 4.7 and 4.8 from~\cite{ARS}, we can
change the mapping $F$ on $Y\setminus A'$ in such a way
that new mapping $F'\colon Y\to l_2$ has the following properties:
\begin{itemize}
\item[(1)] the restriction of $F'$ to the fiber $f^{-1}(x)$
is an embedding for all $x\in X\setminus A$;
\item[(2)] $\dist(F(y),F'(y))< \beta(y)$ for all $y\in Y\setminus A'$.
\end{itemize}

\noindent Upper semicontinuity of $F'$ easily follows from (2).
Let the map $\Psi'$ be given by
$\Psi'(x)=F'\circ f^{-1}(x)$ for all $x\in X\setminus A$.
>From (1) it follows that $\Psi'(x)$ is homeomorphic to
$UV^{[L]}$-compactum $f^{-1}(x)$ for all $x\in X\setminus A$.
Clearly, $\Psi'$ is upper semicontinuous.
\end{proof}

A proper continuous mapping with preimages of points being
$UV^{[L]}$-compacta is called {\it $UV^{[L]}$-mapping}.
The following factorization theorem is known for $n$-soft maps~\cite{B1}.

\begin{thm} \label{thmfac}
Let $L$ be a finite $CW$-complex such that $[L]\le [S^n]$ for some $n$.
If the composition $f\circ g$ of mappings of Polish spaces is (locally)
$[L]$-soft and $g$ is $UV^{[L]}$-map, then $f$ is (locally) $[L]$-soft.
\end{thm}

\begin{proof}
Let $g\colon Y\to E$ and $f\colon E\to X$ be given maps.
Consider a mapping $\varphi\colon Z\to X$ of Polish space $Z$
with $\ed Z\le [L]$ and a mapping $\psi\colon A\to E$ of a
closed subset $A\subset Z$ such that $f\circ \psi=\varphi|_A$.

A multivalued mapping $\Phi=g^{-1}\circ f^{-1}\circ \varphi\colon Z\to Y$
is complete and lower $[L]$-continuous by Theorem~\ref{thmsoftness}.
We have u.s.c. $UV^{[L]}$-valued submapping
$\Psi=g^{-1}\circ \psi\colon A\to Y$ of the map $\Phi$.
By Lemma~\ref{lemmaextUVL} there is u.s.c. $UV^{[L]}$-valued submapping
$\Psi'$ of $\Phi$ defined on some neighbourhood $OA$ of $A$
such that $\Psi'|_A=\Psi$ and $\Psi'|_{OA\setminus A}$ is
continuous and singlevalued.
Clearly, if the map $f\circ g$ is $[L]$-soft, we may assume $OA=Z$.
Then the mapping $\psi'=g\circ\Psi'$ extending $\psi$ is singlevalued
and continuous, and $f\circ \psi'=\varphi|_{OA}$.
\end{proof}

The following corollary was known for $L = S^{k}$
(see \cite[Propositions 2.1.1 and 2.1.2(ii)]{bestvina}).

\begin{cor}
Let $L$ be a finite $CW$-complex such that $[L]\le [S^n]$ for
some $n$ and $g \colon X \to Y$ be a $UV^{[L]}$-map between
Polish spaces. If $X \in A(N)E([L])$, then $Y \in A(N)E([L])$.
\end{cor}
\begin{proof}
Apply Theorem \ref{thmfac} to the composition $f \circ g$,
where $f \colon Y \to \{\text{pt}\}$ is a constant map.
\end{proof}

\begin{thm}\label{T:hurewicz}
Let $f\colon X\to Y$ be a mapping of metric compacta where $\dim Y<\infty$.
Suppose that $\ed Y \leq [M]$ for some finite $CW$-complex~$M$.
If for some locally finite countable $CW$-complex $L$ we have
$\ed (f^{-1}(y)\times Z)\le [L]$ for every point $y\in Y$ and any
Polish space $Z$ with $\ed Z\le [M]$, then $\ed X\le [L]$.
\end{thm}

\begin{proof}
Suppose $A\subset X$ is closed and $g\colon A\to L$ is a map.
We are going to find a continuous extension $\wt g\colon X\to L$ of $g$.
Let $K$ be the cone over $L$ with a vertex $v$.
Denote $\mathcal W=\{h\in C(X,K)\mid h|_A=g\}$ --- a closed
subspace of $C(X,K)$. We define a multivalued map
$F\colon Y\to \mathcal W$ as follows:
$$  F(y)=\{h\in C(X,K)\mid h(f^{-1}(y))\subset K\setminus \{v\}\}.  $$

{\it Claim.} $F$ admits continuous singlevalued selection.

If $\varphi\colon Y\to\mathcal W$ is a continuous selection for $F$, then the
mapping $h\colon X\to K$ defined by $h(x)=\varphi(f(x))(x)$ is continuous.
Since $\varphi(f(x))\in F(f(x))$ for every $x\in X$,
we have $h(X)\subset K\setminus \{v\}$.
Now if $\pi\colon K\setminus \{v\}\to L$ denotes the
natural retraction, then $\wt g =\pi\circ h\colon X\to L$
is the desired continuous extension of $h$.

{\it Proof of the claim.}
Since $K$ is Polish space, the space $C(X,K)$ is also Polish
as well as its closed subspace $\mathcal W$.
Clearly, the graph of $F$ is open in $Y\times \mathcal W$,
therefore $F$ is complete.
Lower $[M]$-continuity of $F$ easily follows from the facts that
the space $\mathcal W$ is locally contractible and $F$ has open graph.

Let us prove that the inclusion $F\subset F$ is fiberwise $[M]$-connected.
Fix a point $y\in Y$ and consider a mapping $\sigma\colon B\to F(y)$
of closed subspace $B$ of a space $Z$ with $\ed Z\le [M]$.
Since $F(y)$ is Polish space, by Corollary~\ref{corL}
we may assume that $Z$ is a Polish space.
It defines a mapping $s\colon B\times X\to K$ by the formula
$s(\{b\}\times \{x\})=\sigma(b)(x)$.
Extend $s$ to a set $Z\times A$ letting $s(\{z\}\times \{a\})=g(a)$.
Clearly, $s$ takes the set $Z\times f^{-1}(y)\cap (Z\times A\cup B\times X)$
into $K\setminus \{v\}\cong L\times [0,1)$.
Since $\ed(Z\times f^{-1}(y))\le [L]$, we can extend $s$ over the set
$Z\times f^{-1}(y)$ to take it into $K\setminus \{v\}$.
Finally extend $s$ over $Z\times X$ as a mapping into $AE$-space $K$.
Now define an extension $\sigma'\colon Z\to F(y)$ of the mapping $\sigma$
by the formula $\sigma'(z)(x)=s(\{z\}\times \{x\})$.

To find a continuous selection of $F$ we apply Theorem~\ref{thmfilteredselection}
to an $n$-filtration $F\subset F\subset\dots\subset F$.
\end{proof}

\appendix
\section{ }

Let $L$ be a $CW$-complex. A pair of spaces $X\subset Y$ is said
to be {\it $[L]$-connected for Polish spaces}
if for every Polish space $Z$ of extension dimension $\ed Z\le [L]$
and for every closed subspace $T \subset Z$ any mapping of $T$
into $X$ can be extended to a mapping of $Z$ into $Y$.

\begin{pro}\label{P:A1}
Let $L$ be a countable locally finite $CW$-complex and $X \subseteq Y$
be a $[L]$-connected pair for Polish spaces in which $X$ is a Polish space.
Then for every completely regular space
$Z$ of extension dimension $\ed Z\le [L]$
and for every $C$-embedded subspace $T\subset Z$ any mapping of $T$
into $X$ can be extended to a mapping of $Z$ into $Y$. In other words,
$X \subseteq Y$ is
$[L]$-connected in the sense of Definition \ref{D:lconnected}.
\end{pro}
\begin{proof}
Consider the Hewitt realcompactification $\upsilon Z$ of the space $Z$. Note that
$\ed \upsilon Z \leq [L]$ (see~\cite{Ch97}, \cite[Theorem 5.1]{Ch}).
By \cite[Theorem 5.2]{Ch}, the realcompact space $\upsilon Z$
can be represented as the limit space of a Polish spectrum
$\displaystyle {\mathcal S}_{\nu Z} = \{ Z_{\alpha}, p_{\alpha}^{\beta}, A\}$ such that
$\ed Z_{\alpha} \leq [L]$ for each $\alpha \in A$.

Since $T$ is $C$-embedded in $Z$ it follows that
$\operatorname{cl}_{\upsilon Z}T$ coincides with the Hewitt
realcompactification $\upsilon T$ of $T$. Next consider the inverse spectrum
$\displaystyle {\mathcal S}^{\prime} =
\{ \operatorname{cl}_{Z_{\alpha}}p_{\alpha}(T),
q_{\alpha}^{\beta}, A\}$, where
$\displaystyle q_{\alpha}^{\beta} = p_{\alpha}^{\beta}\left| \operatorname{cl}_{Z_{\alpha}}p_{\alpha}(T)\right.$
for each $\alpha,\beta \in A$ with $\alpha \leq \beta$.
Since $\upsilon T$ is closed in $\upsilon Z$ it follows that
$\lim{\mathcal S}^{\prime} = \upsilon T$. It is clear that
$\upsilon T$ is $C$-embedded in $\upsilon Z$. This observation,
combined with the fact that the spectrum
${\mathcal S}$ is factorizing, guarantees that the spectrum
${\mathcal S}^{\prime}$ is also factorizing. Now consider a
continuous mapping $f \colon T \to X$. Since $X$ is Polish
there exists a continuous extension $\tilde{f} \colon \upsilon T \to X$.
Factorizability of the spectrum ${\mathcal S}^{\prime}$ implies that we can find an index $\alpha \in A$ and a continuous mapping $f_{\alpha} \colon \operatorname{cl}_{Z_{\alpha}}p_{\alpha}(T) \to X$ such that
$\tilde{f} = f_{\alpha}\circ p_{\alpha}\left| \upsilon T\right.$.
Now recall that the pair $X \subseteq Y$ is $[L]$-connected
and that $Z_{\alpha}$ is a Polish space such that $\ed Z_{\alpha} \leq [L]$.
Consequently there exists a continuous extension $g_{\alpha} \colon Z_{\alpha} \to Y$
of $f_{\alpha}$. Finally consider the composition $p_{\alpha}\circ g_{\alpha} \colon \upsilon Z \to Y$ and let $g = p_{\alpha}\circ g_{\alpha}\left| Z\right.$. Straightforward verification shows that $f = g\left| T\right.$.
\end{proof}

Since every closed subspace of any normal space is $C$-embedded
in it we obtain the following statement.

\begin{cor} \label{corL}
Let $X \subseteq Y$ be a $[L]$-connected pair of Polish spaces.
Then for every paracompact space
$Z$ of extension dimension $\ed Z\le [L]$
and for every closed subspace $T\subset Z$ any mapping of $T$
into $X$ can be extended to a mapping of $Z$ into $Y$.
\end{cor}

The following statement also can be proved using the spectral
technique as presented in \cite{Ch96} (compare to the proof of
Proposition \ref{P:A1}).

\begin{pro}\label{polishane}
Let $L$ be a countable locally finite $CW$-complex
and $X$ be a Polish space.
If $X\in ANE ([L])$ for Polish spaces, then $X\in ANE ([L])$.
\end{pro}


\begin{thebibliography}{99}

\bibitem{ARS} S.~M.~Ageev, D.~Repov\v s, E.V.~\v S\v cepin,
{\it The extension problem for complete $UV\sp n$-preimages},
Tsukuba J. Math. 23 (1999), 97--111.

\bibitem{bestvina}
M.~Bestvina, {\em Characterizing $k$-dimensional universal Menger compacta},
Mem. Amer. Math. Soc. 71 (380), 1988.

\bibitem{Bor} K.~Borsuk,
{\it Theory of Retracts}, PWN, Warszawa, 1967.

\bibitem{B} N.~B.~Brodsky,
{\it On extension of compact-valued mappings},
Uspekhi Mat. Nauk, 54 (1999), 251--252.

\bibitem{B1}  N.~B.~Brodsky
{\it On extension of $UV^n$-valued mappings},
Mat. Zametki, 1999, vol. 66, no. 3, pages 351--363 (in Russian).

\bibitem{B2} N.~B.~Brodsky,
{\it Sections of maps with fibers homeomorphic to a two-dimensional manifold}
Topology Appl., to appear.

\bibitem{BCh} N.~B.~Brodsky, A.~Chigogidze,
{\it Extension dimensional approximation theorem},
Preprint math.GN/0103061 (2001).

\bibitem{Ch97} A.~Chigogidze,
{\it Cohomological dimension of Tychonov spaces},
Topology Appl. 79 (1997), 197-228.

\bibitem{Ch} A.~Chigogidze,
{\it Infinite Dimensional Topology and Shape Theory},
in: Handbook of Geometric Topology, ed. by R.~Daverman and R.~Sher
(2001), to appear.

\bibitem{Ch96} A.~Chigogidze, {\em Inverse Spectra}, North Holland,
Amsterdam, 1996.

\bibitem{ChV} A.~Chigogidze, V.~Valov,
{\it Extension dimension and C-spaces},
Preprint math.AT/0002001 (2001).

\bibitem{Dr1} A.N.~Dranishnikov,
{\it Absolute extensors in dimension $n$ and $n$-soft
mappings increasing the dimension}, Uspekhi Mat. Nauk 39 (1984), 55--95.

\bibitem{Dr2} A.~N.~Dranishnikov,
{\it The Eilenberg-Borsuk theorem for maps in an arbitrary complex},
Russian Acad. Sci. Sb. Math. {\bf 81} (1995), 467--475.

\bibitem{DrD} A.~N.~Dranishnikov, J.~Dydak,
{\it Extension dimension and extension types},
Proc. Steklov Inst. Math. 212 (1996), 55--88.

\bibitem{DRS} A.~N.~Dranishnikov, D.~Repov\v s, E.~V.~\v S\v cepin,
{\it Transversal intersection formula for compacta},
Topology Appl. 85 (1998), 93--117.

\bibitem{DM} J.~Dugundji, E.~Michael,
{\it On local and uniformly local topological properties},
Proc. Amer. Math. Soc. 7 (1956), 304--307.

\bibitem{Eng} R.~Engelking,
{\it Dimension theory}, PWN, Warsaw, 1978.


\bibitem{Kr1} W.~Kryszewski,
{\it Graph-approximation of set-valued maps. A survey},
Differential Inclusions and Optimal Control.
Lecture Notes in Nonlinear Analysis (1998), 223--235.

\bibitem{Le} M.~Levin,
{\it On extensional dimension of maps},
Topology Appl. 103 (2000), 33--35.

\bibitem{LL} M.~Levin, W.~Lewis,
{\it Some mapping theorems for extensional dimension},
Preprint math.GN/0103199 (2001).

\bibitem{Mi89} E.~Michael,
{\it A Generalization of a Theorem on Continuous Selections},
Proc. Amer. Math. Soc., 105 (1989), 236--243.

\bibitem{RS} D.~Repov\v s, P.~V.~Semenov,
{\it Continuous selections of multivalued mappings},
Kluwer Academic Publishers, Dordrecht, 1998.

\bibitem{SB} E.~V.~Shchepin, N.~B.~Brodsky,
{\it Selections of filtered multivalued mappings},
Proc. Steklov Inst. Math., 212 (1996), 209--229.

\end{thebibliography}
\end{document}